%% file: master.tex
\pdfoutput=1
\documentclass{amsart}
\usepackage[utf8]{inputenc}

\usepackage[margin=1in]{geometry}
\usepackage[expansion=false]{microtype}

\usepackage{amsmath, amsfonts, amssymb, amsthm, amsopn}
\usepackage{eucal}

\usepackage{tikz}
\usetikzlibrary{diagrams}

\usepackage[pdfusetitle,unicode]{hyperref}

\newcommand{\resp}{{\sfcode`\.1000 resp.}}
\newcommand{\ie}{{\sfcode`\.1000 i.e.}}
\newcommand{\eg}{{\sfcode`\.1000 e.g.}}

\numberwithin{equation}{section}

\theoremstyle{plain}
\newtheorem*{theorem*}{Theorem}

\newtheorem{theorem}[equation]{Theorem}
\newtheorem{proposition}[equation]{Proposition}
\newtheorem{lemma}[equation]{Lemma}
\newtheorem{corollary}[equation]{Corollary}

\theoremstyle{definition}
\newtheorem{definition}[equation]{Definition}
\newtheorem{example}[equation]{Example}
\newtheorem*{acknowledgements}{Acknowledgements}
\newtheorem*{conventions}{Conventions}

\theoremstyle{remark}
\newtheorem{remark}[equation]{Remark}

\let\scr=\mathcal
\let\bb=\mathbb

\def\Z{\bb Z}
\def\Q{\bb Q}

\def\A{\bb A}
\def\P{\bb P}

\def\1{\mathbf 1}

\let\del=\partial
\let\from=\leftarrow
\let\into=\hookrightarrow
\let\onto=\twoheadrightarrow
\let\tens=\otimes

\def\ph{\mathord-}
\def\minus{-}
\def\brac#1{\langle #1\rangle}
\let\J=\brac
\def\d{\,d}

\def\fp{{\mathrm{fp}}}

\def\qc{{\mathrm{qc}}}

\def\id{\mathrm{id}}
\DeclareMathOperator{\Hom}{Hom}
\DeclareMathOperator{\Fun}{Fun}
\DeclareMathOperator{\End}{End}
\DeclareMathOperator{\Aut}{Aut}
\DeclareMathOperator{\Map}{Map{}}
\DeclareMathOperator{\Spec}{Spec}
\DeclareMathOperator{\Sym}{Sym}
\def\Th{\mathrm{Th}}
\DeclareMathOperator{\GW}{GW}
\DeclareMathOperator{\W}{W}
\DeclareMathOperator{\rk}{rk}
\def\Nis{\mathrm{Nis}}
\def\Zar{\mathrm{Zar}}
\def\et{\mathrm{\acute et}}

\DeclareMathOperator{\Char}{char}
\DeclareMathOperator{\Tr}{Tr}
\DeclareMathOperator{\tr}{tr}

\def\Ex{\mathit{Ex}}
\def\Pr{\mathit{Pr}}
\def\MW{\mathrm{MW}}

\def\PSh{\mathrm{PSh}{}}

\def\Shv{\mathrm{Shv}{}}
\def\D{\mathrm D}
\def\SH{\mathrm{SH}{}}
\def\H{\mathrm H}
\def\Sm{\mathrm{Sm}{}}

\def\op{\mathrm{op}}

\def\ev{\mathrm{ev}}
\def\coev{\mathrm{coev}}

\let\lim=\relax
\DeclareMathOperator*{\lim}{lim}
\DeclareMathOperator*{\colim}{colim}

\title{A quadratic refinement of the Grothendieck--Lefschetz--Verdier trace formula}
\author{Marc Hoyois}
\date{\today}
\address{Department of Mathematics, Massachusetts Institute of Technology, Cambridge, MA, USA}
\email{hoyois@mit.edu}
\urladdr{\url{http://math.mit.edu/~hoyois/}}

\begin{document}

\begin{abstract}
	We prove a trace formula in stable motivic homotopy theory over a general base scheme, equating the trace of an endomorphism of a smooth proper scheme with the ``Euler characteristic integral'' of a certain cohomotopy class over its scheme of fixed points.
	 When the base is a field and the fixed points are étale, we compute this integral in terms of Morel's identification of the ring of endomorphisms of the motivic sphere spectrum with the Grothendieck–Witt ring. In particular, we show that the Euler characteristic of an étale algebra corresponds to the class of its trace form in the Grothendieck–Witt ring.
\end{abstract}

\maketitle
\tableofcontents

\input lefschetz.tex

\newpage

\providecommand{\bysame}{\leavevmode\hbox to3em{\hrulefill}\thinspace}

\end{document}

%% file: lefschetz.tex
\newpage
\section{Introduction and examples}
\label{sec:introduction}

Let $k$ be a field, $X$ a smooth proper $k$-scheme, and $f\colon X\to X$ a $k$-morphism. The Grothendieck–\<Lefschetz–\<Verdier trace formula, originally proved in \cite[Exposé III, \S4]{SGA5}, identifies the trace of the action of $f$ on the $\ell$-adic cohomology of $X$ with the integral of a cohomology class on the scheme of fixed points $X^f$. In the special case where $X^f$ is étale over $k$, the trace formula takes the following simple form:

\begin{theorem}
	\label{thm:classical}
	Let $k$ be a field, $X$ a smooth and proper $k$-scheme, and $f\colon X\to X$ a $k$-morphism with étale fixed points. Then
	\[\sum_{i}(-1)^i\tr (f^\ast\vert H^i_\ell(\bar X))=\sum_{x\in X^f}[\kappa(x):k],\]
	where $\bar X$ is the pullback of $X$ to an algebraic closure of $k$, $\ell\neq\Char k$ is a prime number, and $H^*_\ell(\ph)$ is $\ell$-adic cohomology with coefficients in $\Q_\ell$.
\end{theorem}

The trace formula is thus an equality between two integers associated with $f$. The starting point of the present article is the observation that the left-hand side of the trace formula has a canonical refinement to an element of the \emph{Grothendieck–Witt ring} $\GW(k)$ of the field $k$. To explain why, we need to recall some facts from stable motivic homotopy theory.

Let $\Sm_k$ be the category of smooth separated schemes of finite type over $k$. Consider the functor
\[C_*^\ell\colon \Sm_k\to \widehat\D(\Spec k_\et,\Z_\ell)\]
which sends $p\colon X\to\Spec k$ to the $\ell$-adic sheaf $p_!p^!\Z_\ell$ on $\Spec k_\et$ (here $\widehat\D(B_\et,\Z_\ell)$ is the $\infty$-categorical limit over $n\geq 0$ of the derived categories $\D(B_\et,\Z/\ell^n)$).
By standard properties of $\ell$-adic cohomology and the definition of the stable motivic homotopy category $\SH(k)$,
there is a canonical factorization
\begin{tikzmath}
	\def\colsep{2em}
	\diagram{\Sm_k & \widehat\D(\Spec k_\et,\Z_\ell) \\ \SH(k) & \\};
	\arrows (11) edge node[left]{$\Sigma^\infty_+$} (21) (11-) edge node[above]{$C_*^\ell$} (-12) (21) edge[dashed] node[below right]{$R_\ell$} (12);
\end{tikzmath}
where $R_\ell$ is a symmetric monoidal functor. The functor $\Sigma^\infty_+$ satisfies a generalized version of Poincaré duality, which asserts in particular that, if $X$ is smooth and proper over $k$, $\Sigma^\infty_+X$ is strongly dualizable. Thus, if $f\colon X\to X$ is a $k$-morphism, $\Sigma^\infty_+ f$ has a trace $\tr(\Sigma^\infty_+ f)$ which is an endomorphism of the motivic sphere spectrum $\1_k\in\SH(k)$. Since symmetric monoidal functors commute with traces,
\[R_\ell(\tr(\Sigma^\infty_+ f))=\tr(R_\ell(\Sigma^\infty_+ f))=\tr(C_*^\ell f),\]
and it is clear that $\tr(C_*^\ell f)$ equals the alternating sum appearing in Theorem~\ref{thm:classical}.

Recall that $\GW(k)$ is the group completion of the semiring of isomorphism classes of nondegenerate symmetric bilinear forms over $k$ (or equivalently of nondegenerate quadratic forms if $\Char k\neq 2$).
Associating to such a form the rank of its underlying vector space defines a ring homomorphism
\[\rk\colon\GW(k)\to \Z\]
which is an isomorphism if and only if $k$ is quadratically closed. Given $u\in k^\times$, we denote by $\brac u$ the class of the symmetric bilinear form $k\times k\to k$, $(a,b)\mapsto uab$. These basic classes generate $\GW(k)$ as a group. A fundamental result of Morel\footnote{This result is proved in \cite{Morel} under the assumption that $k$ is perfect. However, Morel actually computes the Nisnevich sheaf on $\Sm_k$ associated with the presheaf $X\mapsto [\Sigma^\infty_+X,\1_k]$, and combining this stronger result with the base change arguments from \cite[Appendix A]{Hoyois} allows us to remove the assumption on $k$.} states that there is a natural isomorphism
\begin{equation}\label{eqn:Morel}
	\GW(k)\simeq\End(\1_k).
\end{equation}
To describe Morel's isomorphism, we first consider a more general construction. Suppose that $V$ is a vector bundle over a scheme $X$ and that $\phi\colon V\stackrel\sim\to V$ is a linear automorphism of $V$. The vector bundle $V$ induces a self-equivalence $\Sigma^V$ of $\SH(X)$, which can be informally described as ``smash product with the sphere bundle of $V$''. The composition
\[\1_X\simeq \Sigma^{-V}\Sigma^{V}\1_X\xrightarrow{\Sigma^{\phi}} \Sigma^{-V}\Sigma^{V}\1_X\simeq \1_X\]
is an automorphism of the motivic sphere spectrum over $X$, which we denote by $\J{\phi}$.\footnote{This construction is of course the algebro-geometric analog of the $J$-homomorphism.}
The isomorphism~\eqref{eqn:Morel} is then given by sending $\brac{u}$ to $\J{u}$, viewing $u\in k^\times$ as a linear automorphism of $\A^1_k$. 
Putting all these facts together, we can identify $\tr(\Sigma^\infty_+ f)$ with a lift of the integer $\tr(C_*^\ell f)$ to $\GW(k)$. It is then natural to ask whether the right-hand side of the Grothendieck–\<Lefschetz–\<Verdier trace formula also lifts to $\GW(k)$, \ie, whether there exist fixed-point indices $i(f,x)\in\GW(k)$, of rank $[\kappa(x):k]$, such that \[\tr(\Sigma^\infty_+f)=\sum_{x\in X^f}i(f,x).\]
An affirmative answer is given in Corollary~\ref{cor:trace} below.\footnote{The existence of such a fixed-point formula was mentioned by Morel in \cite[Remark 4.12 (2)]{MorelICM}.}
It is a consequence of some more general results which we now discuss.

We consider an arbitrary
base scheme $B$. If $X$ is a smooth $B$-scheme such that $\Sigma^\infty_+X\in\SH(B)$ is strongly dualizable, \eg, a smooth proper $B$-scheme, we write
\[\chi(X)=\tr(\Sigma^\infty_+\id_X)\in\End(\1_B)\]
for its Euler characteristic in $\SH(B)$.
More generally, if $\omega$ is an endomorphism of $\1_X$ in $\SH(X)$,
we define
 \[\int_X\omega\d\chi=\tr(p_\sharp \omega)\in\End(\1_B),\]
where $p\colon X\to B$ is the structure map and $p_\sharp\colon\SH(X)\to\SH(B)$ is left adjoint to the base change functor $p^*$.
Note that, by this adjunction, an endomorphism of $\1_X$ is the same thing as a morphism $\Sigma^\infty_+X\to\1_B$ in $\SH(B)$.
The map $\omega\mapsto\int_X\omega \d\chi$ is thus an $\End(\1_B)$-linear functional on the algebra of $\1_B$-valued functions on $\Sigma^\infty_+X$, such that $\int_X1\d\chi = \chi(X)$.
 
We can now state the main result of this paper. Let $X$ be a smooth $B$-scheme, $f\colon X\to X$ a $B$-morphism, and $i\colon X^f\into X$ the inclusion of the scheme of fixed points of $f$.
We say that $f$ has \emph{regular fixed points} if
\begin{itemize}
	\item $X^f$ is smooth over $B$ and
	\item the endomorphism of the conormal sheaf $\scr N_i$ induced by $\id-i^\ast(df)$ is an isomorphism.
\end{itemize}

\begin{theorem}\label{thm:lefschetz}
	Let $X$ be a smooth and proper $B$-scheme and $f\colon X\to X$ a $B$-morphism with regular fixed points. Then \[\tr(\Sigma^\infty_+f)=\int_{X^f}\J{\phi}\d\chi,\]
	where $\phi$ is the automorphism of the conormal sheaf of the immersion $i\colon X^f\into X$ induced by $\id-i^*(df)$.
\end{theorem}

Theorem~\ref{thm:lefschetz} will be proved in~\S\ref{sec:traces}.
The following special case is worth recording:

\begin{corollary}\label{cor:lefschetz}
	Let $X$ be a smooth and proper $B$-scheme and $f\colon X\to X$ a $B$-morphism. If $\tr(\Sigma^\infty_+f)\neq 0$, then $f$ has a fixed point.
\end{corollary}

Along the way we will observe that $\chi(X)=0$ if $\Omega_{X/B}$ has a nonvanishing global section (see Remark~\ref{rmk:vanishing}):

\begin{theorem}\label{thm:vanishing}
	Let $X$ be a smooth and proper $B$-scheme. If $[\Omega_{X/B}]=[\scr O_X]+[\scr E]$ in $K_0(X)$ for some locally free sheaf $\scr E$, then $\int_X\omega\d\chi=0$ for all $\omega\in\End(\1_X)$.
\end{theorem}

The properness hypothesis in Theorems~\ref{thm:lefschetz} and~\ref{thm:vanishing} is essential: there are many smooth $B$-schemes that become strongly dualizable in $\SH(B)$ without being proper (\eg, the complement of a smooth closed subscheme in a smooth proper scheme), but these theorems clearly do not extend to all such schemes.

Before giving examples, we make some general remarks on the notion of regular fixed points appearing in Theorem~\ref{thm:lefschetz}. Let $\Delta_X\subset X\times_BX$ be the diagonal and let $\Gamma_f\subset X\times_BX$ be the graph of $f$. It is clear that we have the following implications:
\[\Gamma_f\text{ and }\Delta_X\text{ intersect transversely}
\:\Rightarrow\:
f\text{ has regular fixed points}
\:\Rightarrow\:
\Gamma_f\text{ and }\Delta_X\text{ intersect cleanly}\]
(the latter simply means that $X^f$ is smooth over $B$). Moreover, both implications are strict: the transposition on $X\times_BX$ has regular fixed points if and only if multiplication by $2$ on $\Omega_{X/B}$ is invertible. Even in the case of a transverse intersection, we will see in Example~\ref{ex:euler} below that $\int_{X^f}\J\phi\d\chi$ can depend on $\phi$. In particular, the trace of $\Sigma^\infty_+f$ is not determined by the derived fixed points of $f$, since the latter coincide with the underived fixed points when the intersection of $\Gamma_f$ and $\Delta_X$ is transverse. This is a significant difference between stable motivic homotopy and $\ell$-adic cohomology. 
\begin{example}[Fixed points of Frobenius]
Let $q$ be a prime power, $X$ a smooth and proper $\bb F_q$-scheme, and $f\colon X\to X$ the Frobenius endomorphism. Then
	\[X^f\simeq\coprod_{X(\bb F_q)}\Spec\bb F_q\]
	and $df=0$. By Theorem~\ref{thm:lefschetz}, $\tr(\Sigma^\infty_+f)\in\GW(\bb F_q)$ is simply the Euler characteristic of $X^f$, which is the number of $\bb F_q$-rational points of $X$ by additivity of the trace.
\end{example}

\begin{example}[The Euler characteristic of $\P^1$]\label{ex:euler}
	We can compute the Euler characteristic of projective space $\P^n$ by induction on $n$ using the cofiber sequence
		\[\Sigma^\infty_+\P^{n-1}\to \Sigma^\infty_+\P^n\to S^{\A^n}\]
		and the additivity of the trace (see \cite{May}). We find that
		\[\chi(\P^n)=\begin{cases}
			\chi(\P^{n-1})+1 & \text{if $n$ is even,}\\
			\chi(\P^{n-1})+\tau & \text{if $n$ is odd,}
		\end{cases}\]
		where $\tau\in\End(\1_B)$ is the desuspension of the transposition $S^{\A^1}\wedge S^{\A^1}\simeq S^{\A^1}\wedge S^{\A^1}$.\footnote{Here we use the following fact: if $\scr C$ is a symmetric monoidal category and $L\in\scr C$ is $\tens$-invertible, then $\chi(L)\in\End(\1)$ corresponds to the transposition under the canonical isomorphism $\End(L\tens L)\simeq \End(\1)$. We leave the elementary proof to the reader.} If $B$ is the spectrum of a field $k$, it is well-known that $\tau$ corresponds to $\brac{-1}\in\GW(k)$.
	As a consistency test, we use Theorem~\ref{thm:lefschetz} to show that the Euler characteristic of the projective line $\P^1$ over $k$ is the hyperbolic form $\brac{1, -1}\in\GW(k)$. Since an odd-degree extension of finite fields induces an isomorphism on Grothendieck–Witt rings, we may assume without loss of generality that $k$ has at least $4$ elements. Choose $a\in k^\times$ with $a^2\neq 1$ and let $f$ be the automorphism of $\P^1$ given by $[x:y]\mapsto[a^2x:y]$. A homotopy between the matrices
	\[\begin{pmatrix}1 & 0 \\ 0 & 1 \end{pmatrix}\quad\text{and}\quad \begin{pmatrix}a & 0 \\ 0 & a^{-1} \end{pmatrix}\]
	in $\mathrm{SL}_2(k)$ induces a homotopy between $\id_{\P^1}$ and $f$, so that $\chi(\P^1)=\tr(\Sigma^\infty_+f)$. We have
	\[(\P^1)^f=\{0,\infty\}\]
	(a disjoint union of two copies of $\Spec(k)$), $df_0=a^2$, and $df_\infty=a^{-2}$. Thus, the endomorphism $\id-i^\ast(df)$ of $i^\ast(\Omega_{\P^1})\simeq \scr N_i$ is multiplication by $1-a^2$ at $0$ and by $1-a^{-2}$ at $\infty$. By Theorem~\ref{thm:lefschetz}, the trace of $f$ is $\brac{1-a^2, 1-a^{-2}}=\brac{1, -1}$, as expected.
\end{example}

\begin{example}[Relations in the endomorphism ring of the motivic sphere spectrum]
	The fact that $\tr(\Sigma^\infty_+f)$ is an invariant of the homotopy class of $f$ produces interesting relations in the ring $\End(\1_B)$. For example, if $k$ is a field and $a_0,\dotsc,a_n\in k^\times$ are $n+1$ distinct elements whose product is $1$, then the endomorphism $[x_0:\dotsc:x_n]\mapsto[a_0x_0:\dotsc:a_nx_n]$ of $\P^n$ over $k$ is homotopic to the identity. It follows that its trace, which by Theorem~\ref{thm:lefschetz} is the class
	\[\sum_{i=0}^n\prod_{j\neq i}\brac{1-a_j/a_i}\in\GW(k),\]
	is independent of the choice of the elements $a_i$ and equals the Euler characteristic of $\P^n$.
\end{example}

Our proofs of Theorems~\ref{thm:lefschetz} and~\ref{thm:vanishing} remain valid if the functor $B\mapsto\SH(B)$ is replaced by any \emph{motivic triangulated category} in the sense of \cite[Definition 2.4.45]{CD}. On the other hand, by the $\infty$-categorical universality of $\SH(B)$ for \emph{fixed} $B$ established in \cite{Robalo}, our theorems admit the following generalizations. Let $\scr C$ be a pointed symmetric monoidal presentable $\infty$-category and $F\colon \Sm_B\to\scr C$ a symmetric monoidal functor satisfying $\A^1$-homotopy invariance,
Nisnevish descent,
and $\P^1$-stability (\ie, the cofiber of $F(\infty)\to F(\P^1_B)$ is $\tens$-invertible). Then $F$ sends smooth proper $B$-schemes to strongly dualizable objects and Theorems~\ref{thm:lefschetz} and~\ref{thm:vanishing} are true with $\Sigma^\infty_+$ replaced by $F$.
For example, when $B$ is a field and $F=C_*^\ell$, Theorem~\ref{thm:lefschetz} recovers Theorem~\ref{thm:classical}. 
Finally, in \S\ref{sec:fields}, we will prove:

\begin{theorem}\label{thm:trace}
	Let $k\subset L$ be a finite separable field extension, $V$ a finite-dimensional vector space over $L$, and $\phi$ an automorphism of $V$. Then, modulo the isomorphism~\eqref{eqn:Morel},
	\[\int_{L}\J\phi\d\chi=\Tr_{L/k}\brac{\det(\phi)}.\]
\end{theorem}

Here, $\Tr_{L/k}\colon\GW(L)\to\GW(k)$ is the Scharlau transfer associated with the field trace $\Tr_{L/k}\colon L\to k$, \ie, it sends a symmetric bilinear form $b\colon V\times V\to L$ to the form ${\Tr_{L/k}}\circ b\colon V\times V\to k$ of rank $[L:k]\rk(b)$. Note that we allow $k$ to have characteristic $2$ or to be imperfect.
Combining Theorems~\ref{thm:lefschetz} and~\ref{thm:trace} gives the following result, which is a motivic version of the Lefschetz–Hopf theorem \cite[VII, Proposition 6.6]{Dold}: 
\begin{corollary}\label{cor:trace}
Let $k$ be a field, $X$ a smooth and proper $k$-scheme, and $f\colon X\to X$ a $k$-morphism with étale fixed points. Then
\[\tr(\Sigma^\infty_+f)=\sum_{x\in X^f}\Tr_{\kappa(x)/k}\brac{\det(\id-df_x)}.\]
\end{corollary}

\begin{example}[The Euler characteristic of $\P^1$, continued]
	Let $k$ be a field such that $\sqrt{-1}\notin k$. Consider the endomorphism $f\colon\P^1\to\P^1$ given by $[x:y]\mapsto [-y:x]$. It is again induced by a matrix in $\mathrm{SL}_2(k)$ and hence, as in Example~\ref{ex:euler}, is homotopic to $\id_{\P^1}$. We have
	\[(\P^1)^f\simeq \Spec k(i),\]
where $i$ is a square root of $-1$. Moreover, $df_i$ is multiplication by $i^{-2}=-1$. The fixed-point index of $f$ at $i$ is therefore
	\[\Tr_{k(i)/k}\brac{1-(-1)}=\brac{4,-4}=\brac{1,-1}\in\GW(k).\]
	As predicted by Corollary~\ref{cor:trace}, this coincides with the Euler characteristic of $\P^1$ computed in Example~\ref{ex:euler}.
\end{example}

\begin{conventions}
The following conventions are in force throughout, except in Appendix~\ref{app:qcqs}:
\begin{itemize}
	\item All schemes are assumed to be coherent, \ie, quasi-compact and quasi-separated.
	\item Smooth and étale morphisms are assumed to be separated and of finite type.
\end{itemize}
See however Remark \ref{rmk:gluing}.
\end{conventions}

\begin{acknowledgements}
	I thank Marc Levine and Jean Fasel for their interest in this project and for stimulating conversations about it.
	The first version of this paper was written while I was visiting the department of mathematics at the University of Duisburg–Essen and I would like to thank everyone there for their hospitality.
	Finally, I am immensely grateful to the anonymous referee whose report lead to considerable improvements to the original manuscript.
\end{acknowledgements}

\newpage
\section{Review of the formalism of six operations}
\label{sec:review}

To prove Theorem~\ref{thm:lefschetz}, we will use the formalism of six operations ($f^*$, $f_*$, $f_!$, $f^!$, $\wedge$, and $\Hom$) in stable motivic homotopy theory developed by Ayoub in \cite{Ayoub} and revisited by Cisinski and Déglise in \cite{CD}. In this section we briefly review the main features of this formalism, and we introduce several pieces of notation that will be used throughout this paper.

\begin{remark}
	We do not insist that schemes be noetherian and of finite Krull dimension. We explain in Appendix~\ref{app:qcqs} how to extend motivic homotopy theory and the formalism of six operations to arbitrary schemes.
	 \end{remark}

For $B$ a scheme, we denote by $\SH(B)$ the closed symmetric monoidal triangulated category of motivic spectra parametrized by $B$. The monoidal unit, monoidal product, monoidal symmetry, and internal hom in $\SH(B)$ will be denoted by $\1_B$, $\wedge$, $\tau$ and $\Hom$, respectively. We first give a description of the six operations which is independent of the specifics of the category $\SH(B)$.

To any morphism of schemes $f\colon Y\to X$ is associated an adjunction
\[f^*:\SH(X)\rightleftarrows \SH(Y): f_*\]
where $f^*$ is symmetric monoidal.
If $f$ is smooth, $f^*$ also admits a left adjoint denoted by $f_\sharp$.
If $f$ is separated of finite type, there is an exceptional adjunction
\[f_!:\SH(Y)\rightleftarrows\SH(X):f^!\]
and a natural transformation $f_!\to f_*$ which is an isomorphism when $f$ is proper. Each of the assignments $f\mapsto f^*,f_*,f_!,f^!,f_\sharp$ is part of a $2$-functor on the category of schemes. In particular, every commutative triangle of schemes gives rise to various \emph{connection isomorphisms}, such as $(gf)^*\simeq f^*g^*$, satisfying cocyle conditions. We will denote by $c$ any isomorphism which is a composition of such connection isomorphisms.

To any cartesian square of schemes
\begin{tikzequation}\label{eqn:cartesian}
	\diagram{\bullet & \bullet \\ \bullet & \bullet \\};
	\arrows (11-) edge node[above]{$g$} (-12) (11) edge node[left]{$q$} (21) (21-) edge node[below]{$f$} (-22) (12) edge node[right]{$p$} (22);
\end{tikzequation}
are associated several \emph{exchange transformations} such as
\begin{gather*}
	\Ex^*_*\colon f^*p_*\to q_*g^*, \\
	\Ex^{*!}\colon g^*p^!\to q^!f^*, \\
	\Ex_\sharp^*\colon g_\sharp q^*\to p^* f_\sharp.
\end{gather*}
To a morphism $f$ are also associated several \emph{projectors} such as
\begin{gather*}
	\Pr_*^*\colon f_*E \wedge F\to f_*(E\wedge f^* F),\\
	\Pr^{*!}\colon f^*E\wedge f^!F\to f^!(E\wedge F),\\
	\Pr_\sharp^*\colon f_\sharp(E\wedge f^* F)\to f_\sharp E\wedge F.
\end{gather*}
Each projector comes in left and right variants (for which we use the same symbol) related to one another via the monoidal symmetry $\tau$. There are also projectors involving the internal hom, but we will not need them. A crucial fact is that the transformations $\Ex_!^*$, $\Ex^!_*$, and $\Pr_!^*$ are always isomorphisms. As we will see below, this generalizes the proper base change theorem ($\Ex_*^*$ is an isomorphism when $p$ is proper), the smooth base change theorem ($\Ex_*^*$ is an isomorphism when $f$ is smooth), and the projection formula ($\Pr_*^*$ is an isomorphism when $f$ is proper).

If $i\colon Z\into X$ is a closed immersion with open complement $j\colon U\into X$, we have two localization cofiber sequences
\begin{gather*}
	j_!j^!\stackrel\epsilon\longrightarrow \id \stackrel\eta\longrightarrow i_*i^*,\\
	i_!i^!\stackrel\epsilon\longrightarrow\id\stackrel\eta\longrightarrow j_*j^*.
\end{gather*}
Moreover, the functors $i_*\simeq i_!$, $j_!$, and $j_*$ are fully faithful. We will denote by
$\sigma\colon i^!\to i^*$
the natural transformation
\[i^!\simeq \id^*i^!\xrightarrow{\Ex^{*!}} \id^! i^*\simeq i^*.\]

If $p\colon V\to X$ is a vector bundle with zero section $s$, the adjunction
\[p_\sharp s_* : \SH(X)\rightleftarrows \SH(X): s^!p^*\]
is a self-equivalence of $\SH(X)$, which we will denote by $\Sigma^V\dashv \Sigma^{-V}$. The functors $\Sigma^V$ and $\Sigma^{-V}$ will be called \emph{Thom transformations}, or the \emph{$V$-suspension} and \emph{$V$-desuspension} functors, respectively. They are compatible with each of the operations $f^*$, $f_*$, $f_\sharp$, $f_!$, and $f^!$ in the following sense: there are canonical isomorphisms $f^*\Sigma^V\simeq \Sigma^{f^*V}f^*$, $\Sigma^V f_*\simeq f_*\Sigma^{f^*V}$, etc. They are also compatible with the monoidal structure, in the sense that $\Sigma^{V}E\wedge F\simeq\Sigma^V(E\wedge F)$ and $\Sigma^{-V}E\wedge F\simeq\Sigma^{-V}(E\wedge F)$. In particular,
\[\Sigma^V\simeq \Sigma^V\1_X\wedge(\ph)\quad\text{and}\quad\Sigma^{-V}\simeq \Sigma^{-V}\1_X\wedge(\ph).\]
 If $\scr M$ is a locally free sheaf of finite rank on $X$, we will also denote by $\Sigma^{\scr M}$ and $\Sigma^{-\scr M}$ the functors $\Sigma^{\bb V(\scr M)}$ and $\Sigma^{-\bb V(\scr M)}$, where $\bb V(\scr M)=\Spec(\Sym(\scr M))$ is the vector bundle on $X$ whose sheaf of sections is dual to $\scr M$.

If $f$ is smooth, there are canonical isomorphisms
\[f_!\simeq f_\sharp\Sigma^{-\Omega_f}\quad\text{and}\quad f^!\simeq\Sigma^{\Omega_f}f^*,\]
where $\Omega_f$ is the sheaf of relative differentials of $f$. In particular, if $f$ is étale, $f_!\simeq f_\sharp$ and $f^!\simeq f^*$. At this point we see that the operations $f_\sharp$, $\Sigma^V$, and $\Sigma^{-V}$, which are not listed among \emph{the} six operations, are expressible in terms of the latter as follows:
\[f_\sharp\simeq f_!\Sigma^{\Omega_f},\quad \Sigma^V\simeq s^*p^!,\quad \Sigma^{-V}\simeq s^!p^*.\]

The Thom transformations are functorial in monomorphisms of vector bundles (\ie, epimorphisms of locally free sheaves) as follows. Given a triangle
	\begin{tikzmath}
		\def\colsep{.75em}
		\diagram{W & & V \\ & X\rlap, & \\};
		\arrows (11-) edge[c->] node[above]{$\phi$} (-13) (11) edge node[left]{$q$} (22) (13) edge node[right]{$p$} (22);
	\end{tikzmath}
where $p$ and $q$ are vector bundles with zero sections $s$ and $t$ and where $\phi$ exhibits $W$ as a subbundle of $V$, we define $\Sigma^{\phi}\colon \Sigma^{W}\to \Sigma^{V}$ to be the composition
\[t^*q^!\stackrel c\simeq t^*\phi^!p^!\stackrel \sigma\to t^*\phi^*p^!\stackrel c\simeq s^*p^!,\]
and we let $\Sigma^{-\phi}\colon \Sigma^{-V}\to\Sigma^{-W}$ be its mate, which is given by the same composition with stars and shrieks exchanged. In particular, a linear automorphism $\phi\colon V\stackrel\sim\to V$ induces an automorphism $\Sigma^{-V}\Sigma^\phi$ of the identity functor on $\SH(X)$, which we denote by $\J\phi$.

For any short exact sequence
\[0\to W\to V\to U\to 0\]
of vector bundles on $X$, the exchange transformation $\Ex^{*!}$ provides an isomorphism
\begin{equation}\label{eqn:thomses}
	\Sigma^{V}\simeq \Sigma^{W}\Sigma^{U}
\end{equation}
which is natural with respect to monomorphisms of short exact sequence. The properties of these isomorphisms established in \cite[\S1.5]{Ayoub} show that the association $V\mapsto \Sigma^V$ induces a morphism of Picard groupoids
\[\Sigma^{(\ph)}\colon K(X)\to \Aut(\SH(X))\]
from the $K$-theory groupoid of $X$ to the groupoid of self-equivalences of $\SH(X)$. In particular, the map $\phi\mapsto\J\phi$ factors through a group homomorphism $K_1(X)\to \Aut(\id_{\SH(X)})$.

Given a commutative triangle
\begin{tikzequation}\label{eqn:puritytriangle}
	\diagram{Z & X \\ & B\rlap, \\};
	\arrows (11-) edge[c->] node[above]{$s$} (-12) (11) edge node[below left]{$q$} (22) (12) edge node[right]{$p$} (22);
\end{tikzequation}
where $p$ and $q$ are smooth and $s$ is a closed immersion, we obtain a sequence of isomorphisms
\[s^!p^*\simeq s^!\Sigma^{-\Omega_p}p^!\simeq \Sigma^{-s^*(\Omega_p)}s^!p^!\stackrel c\simeq \Sigma^{-s^*(\Omega_p)}q^!\simeq \Sigma^{-s^*(\Omega_p)}\Sigma^{\Omega_q}q^*\simeq \Sigma^{-\scr N_s}q^*,\]
where the last isomorphism is induced by the short exact sequence
\[0\to \scr N_s\to s^*(\Omega_p)\xrightarrow{ds}\Omega_q\to 0.\]
The isomorphism $s^!p^*\simeq \Sigma^{-\scr N_s}q^*$ and its mate $p_\sharp s_*\simeq q_\sharp\Sigma^{\scr N_s}$ are called the \emph{purity isomorphisms} and are denoted by $\Pi$. Although the purity isomorphism appears \textit{a posteriori} as a consequence of the formalism of six operations, it must be constructed ``by hand'' in both the approach of Ayoub and that of Cisinski–Déglise. We discuss the purity isomorphism further in Appendix~\ref{app:purity} (where in particular we show that the constructions of Ayoub and of Ciskinski–Déglise are equivalent).

Of course, all this data satisfies many coherence properties, of which an exhaustive list cannot easily be written down. Let us mention here one kind of coherence that we will use often. If $f$ is a smooth morphism (\resp{} a proper morphism), then we may want to replace, in a given expression, occurrences of $f_!$ and $f^!$ by $f_\sharp\Sigma^{-\Omega_f}$ and $\Sigma^{\Omega_f}f^*$ (\resp{} occurences of $f_!$ by $f_*$). Such replacements yield canonically isomorphic expressions, and, under these canonical isomorphisms, any exchange transformation is transformed into another exchange transformation, and any projector is transformed into another projector. For example, consider the cartesian square~\eqref{eqn:cartesian} and the exchange isomorphism $\Ex_!^*\colon p^* f_!\stackrel\sim\to g_!q^*$. If $f$ is smooth, then $q^*(\Omega_f)\simeq\Omega_g$ and the square
\begin{tikzmath}
	\diagram{p^*f_! & g_!p^* \\ p^*f_\sharp \Sigma^{-\Omega_f} & g_\sharp q^*\Sigma^{-\Omega_f} \\};
	\arrows (11-) edge node[above]{$\Ex_!^*$} (-12) (11) edge node[left]{$\simeq$} (21) (21-) edge[<-] node[below]{$\Ex_\sharp^*$} (-22) (12) edge node[right]{$\simeq$} (22);
\end{tikzmath}
commutes, while the square
\begin{tikzmath}
	\diagram{p^*f_! & g_!p^* \\ p^*f_* & g_* q^* \\};
	\arrows (11-) edge node[above]{$\Ex_!^*$} (-12) (11) edge (21) (21-) edge node[below]{$\Ex_*^*$} (-22) (12) edge (22);
\end{tikzmath}
commutes for any $f$ (the vertical maps being isomorphisms when $f$ is proper). Similarly, when $f$ or $p$ is smooth, the exchange transformation $\Ex^{*!}$ transforms into the isomorphism $\Ex^{!!}$ or $\Ex^{**}$.

Let us now describe these functors more explicitly. For $B$ a scheme, we denote by $\Sm_B$ the category of smooth $B$-schemes and by $\H_{(*)}(B)$ the (pointed) motivic homotopy category of $B$ (we refer to Appendix~\ref{app:qcqs} for the definitions in the generality considered here).
We denote by
\begin{gather*}
	\Sigma^\infty_+\colon \H(B)\to \SH(B),\\
	\Sigma^\infty\colon \H_*(B)\to\SH(B)
\end{gather*}
	the canonical symmetric monoidal functors,
called \emph{stabilization functors}. If $X\in\Sm_B$ and $U\into X$ is an open subscheme, $X/U$ is the quotient of the presheaf represented by $X$ by the presheaf represented by $U$, viewed as an object of $\H_*(B)$. If $V$ is a vector bundle on $X\in\Sm_B$, we denote its \emph{Thom space} by
\[\Th_X(V)=\frac{V}{V\minus X}\in\H_*(B).\]
If $V$ is a vector bundle over $B$ itself, we also write $S^V$ for $\Th_B(V)$ or for its stabilization $\Sigma^\infty\Th_B(V)$.

For $f\colon Y\to X$, the functor $f^*\colon \SH(X)\to\SH(Y)$ is induced by the base change functor $\Sm_X\to \Sm_Y$, so that
\[f^*\Sigma^\infty_+ U\simeq \Sigma^\infty_+(U\times_XY).\]
If $f$ is smooth, the functor $f_\sharp$ is similarly induced by the forgetful functor $\Sm_Y\to\Sm_X$. In particular, if $p\colon X\to B$ is smooth, then 
\[\Sigma^\infty_+X\simeq p_\sharp p^*\1_B\simeq p_!p^!\1_B\in\SH(B).\]
If $i\colon Z\into B$ is a closed immersion with open complement $j\colon U\into B$ and if $X\in\Sm_B$, the localization cofiber sequence
\[j_\sharp \Sigma^\infty_+ X_U\to \Sigma^\infty_+X\to i_* \Sigma^\infty_+X_Z\]
shows that
\[i_*\Sigma^\infty_+X_Z\simeq\Sigma^\infty(X/X_U).\]
In particular, if $V$ is a vector bundle on $X$, then $\Sigma^V\1_X\simeq S^V$ and hence $\Sigma^V\simeq S^V\wedge(\ph)$. If $p\colon X\to B$ is smooth and $V$ is a vector bundle on $X$, we deduce that
\[\Sigma^\infty\Th_X(V)\simeq p_\sharp \Sigma^Vp^*\1_B\simeq p_!\Sigma^Vp^!\1_B\in\SH(B).\]

Consider a commutative triangle
\begin{tikzmath}
	\def\colsep{.75em}
	\diagram{Y & & X \\ & B & \\};
	\arrows (11-) edge node[above]{$f$} (-13) (11) edge node[left]{$q$} (22) (13) edge node[right]{$p$} (22);
\end{tikzmath}
where $p$ and $q$ are smooth. Under the isomorphisms $\Sigma^\infty_+ X\simeq p_!p^!\1_B$ and $\Sigma^\infty_+Y\simeq q_!q^!\1_B$, the map $\Sigma^\infty_+f$ in $\SH(B)$ is given by the composition
\[q_!q^!\1_B\stackrel c\simeq p_!f_!f^!p^!\1_B\stackrel\epsilon\to p_!p^!\1_B\]
(this is \cite[Lemme C.2]{Ayoub2}). More generally, suppose that $V$ and $W$ are vector bundles on $X$ and $Y$ and let $\phi\colon W\into f^*V$ be a monomorphism of vector bundles. Then the map of Thom spectra $\Sigma^\infty\Th_Y(W)\to\Sigma^\infty\Th_X(V)$ induced by $\phi$ is given by the composition
\[q_!\Sigma^Wq^!\1_B\xrightarrow{\Sigma^\phi} q_!\Sigma^{f^*V}q^!\1_B\stackrel c\simeq p_!f_!\Sigma^{f^*V}f^!p^!\1_B\simeq p_!\Sigma^{V}f_!f^!p^!\1_B\stackrel\epsilon\to p_!\Sigma^Vp^!\1_B.\]
This is easily proved by considering the localization cofiber sequences defining $\Th_Y(W)$ and $\Th_X(V)$ and applying the previous result to the maps $W\minus Y\to V\minus X$ and $W\to V$.

Finally, given the triangle~\eqref{eqn:puritytriangle} with $p$ and $q$ smooth and $s$ a closed immersion, the purity isomorphism $\Pi\colon p_\sharp s_*\1_X\simeq q_\sharp\Sigma^{\scr N_s}\1_X$ is the stabilization of the unstable isomorphism
\[\frac{X}{X\minus Z}\simeq\Th_Z(\bb V(\scr N_s))\]
in $\H_*(B)$ from \cite[Theorem 2.23]{MV}.

\newpage
\section{Duality in stable motivic homotopy theory}
\label{sec:duality}

Fix a base scheme $B$. In \cite[Appendix A]{Hu} and \cite[\S2]{Riou} it was proved that smooth and projective
$B$-schemes become strongly dualizable in $\SH(B)$. We will follow the latter reference and deduce this duality as an easy consequence of the formalism of six operations.
We will then provide alternative descriptions of this duality that we will need in \S\ref{sec:traces} and \S\ref{sec:fields}.

Recall that an object $A$ in a symmetric monoidal category $(\scr C,{\tens},\1)$ is \emph{strongly dualizable} if there exists an object $A^\vee$ and morphisms
\[\coev\colon \1\to A\tens A^\vee\quad\text{and}\quad \ev\colon A^\vee\tens A\to\1\]
such that both compositions
\begin{gather*}
	\begin{tikzpicture}[ampersand replacement=\&]
		\def\colsep{4em}
		\diagram{A\simeq \1\tens A \& A\tens A^\vee\tens A \& A\tens \1\simeq A \\};
		\arrows (11-) edge node[above]{$\coev\tens\id$} (-12) (12-) edge node[above]{$\id\tens\ev$} (-13);
	\end{tikzpicture} \\
	\begin{tikzpicture}[ampersand replacement=\&]
		\def\colsep{4em}
		\diagram{A^\vee\simeq A^\vee\tens\1 \& A^\vee\tens A\tens A^\vee \& \1\tens A^\vee\simeq A^\vee \\};
		\arrows (11-) edge node[above]{$\id\tens\coev$} (-12) (12-) edge node[above]{$\ev\tens\id$} (-13);
	\end{tikzpicture}
\end{gather*}
are the identity. When it exists, this data is unique up to a unique isomorphism. If objects $A$ and $A^\vee$ are given, then a choice of coevaluation and evaluation maps exhibiting $A^\vee$ as a strong dual of $A$ is equivalent to a choice of adjunction between $A^\vee\tens(\ph)$ and $A\tens(\ph)$. The counit and unit of such an adjunction determine the evaluation and the coevaluation, respectively. If $A\in\scr C$ is strongly dualizable and $f\colon A\to A$ is an endomorphism, then the \emph{trace} of $f$ is the endomorphism of the unit $\1$ given by the composition
\[\1\xrightarrow{\coev} A\tens A^\vee \xrightarrow{f\tens\id} A\tens A^\vee\stackrel\tau\simeq A^\vee\tens A\xrightarrow{\ev} \1.\]

Throughout this section we fix a smooth and proper morphism $p\colon X\to B$.
Recall that the projector
\[\Pr_!^*\colon p_!(E\wedge p^*F)\to p_!E\wedge F\]
is always an isomorphism. In particular, for $E=p^*\1_B$, we obtain a natural isomorphism
\begin{equation}\label{eqn:duality1}
	p_!p^*\simeq p_!p^*\1_B\wedge(\ph).
\end{equation}
The projectors
\[p_*p^! E\wedge F \xrightarrow{\Pr_*^*}p_*(p^!E\wedge p^* F)\xrightarrow{\Pr^{*!}} p_*p^!(E\wedge F)\]
are also isomorphisms, the first because $p$ is proper and the second because $p$ is smooth. For $E=\1_B$ we obtain an isomorphism
\begin{equation}\label{eqn:duality2}
	p_*p^!\simeq p_*p^!\1_B\wedge (\ph).
\end{equation}
Since $p_!p^*$ is left adjoint to $p_*p^!$, we obtain from \eqref{eqn:duality1} and \eqref{eqn:duality2} an adjunction between $p_!p^*\1_B\wedge(\ph)$ and $p_*p^!\1_B\wedge (\ph)$, \ie, a strong duality between $p_!p^*\1_B$ and $p_*p^!\1_B\simeq \Sigma^\infty_+ X$. Under the isomorphisms \eqref{eqn:duality1} and~\eqref{eqn:duality2}, the coevaluation map $\1_B\to p_*p^!\1_B\wedge p_! p^\ast\1_B$ is the composition of the units
\begin{equation}\label{eqn:coev}
	\1_B\stackrel{\eta}{\to} p_\ast p^\ast\1_B\stackrel{\eta}{\to} p_\ast p^! p_! p^\ast \1_B,
\end{equation}
and the evaluation map $p_! p^\ast\1_B\wedge p_*p^!\1_B\to\1_B$ is the composition of the counits
\begin{equation}\label{eqn:ev}
	p_!p^\ast p_\ast p^!\1_B\stackrel{\epsilon}{\to} p_! p^!\1_B\stackrel{\epsilon}{\to}\1_B.
\end{equation}

\begin{remark}
	Composing the coevaluation with the symmetry and the first half of~\eqref{eqn:ev}, we obtain a map $\1_B\to p_!p^!\1_B\simeq \Sigma^\infty_+ X$ in $\SH(B)$. This is the motivic analog of the Becker–Gottlieb transfer in stable parametrized homotopy theory. 
	It is easy to see that integration against the Euler characteristic is equivalent to precomposition with this transfer.
	\end{remark}

Consider the cartesian square
\begin{tikzmath}
	\diagram{X\times_BX & X \\ X & B\rlap, \\};
	\arrows (11-) edge node[above]{$\pi_2$} (-12) (11) edge node[left]{$\pi_1$} (21)
	(21-) edge node[below]{$p$} (-22) (12) edge node[right]{$p$} (22);
\end{tikzmath}
and denote by $\delta\colon X\into X\times_BX$ the diagonal immersion.
The key result which will be the basis for the proof of the main theorem in \S\ref{sec:traces} is the following description of the trace of an endomorphism:

\begin{proposition}\label{prop:trace}
	Let $f\colon X\to X$ be a $B$-morphism and let $\omega\in\End(\1_X)$. Then $\tr(p_\sharp\omega\circ\Sigma^\infty_+ f)\colon \1_B\to \1_B$ is given by the following composition evaluated at $\1_B$:
	\begin{tikzmath}
		\def\colsep{1.5em}
		\diagram{\id & p_* p^* & p_* \pi_{1!}\delta_!\delta^!\pi_2^! p^* & p_* \pi_{1!}\pi_2^! p^* & & & \\
		& & & p_!\pi_{2*}\pi_1^*p^! & p_!\pi_{2*}\delta_*\delta^*\pi_1^*p^! & p_!p^! & \id\rlap, \\};
		\arrows (11-) edge node[above]{$\eta$} (-12) (12-) edge node[above]{$\simeq$} (-13) (13-) edge node[above]{$\epsilon$} (-14)
		(14) edge[<-] node[left]{$\Ex_{!*}^{\vphantom{!}}$} node[right]{$\Ex^{*!}_{\vphantom{!}}$} (24)
		(24-) edge node[above]{$\eta$} (-25) (25-) edge node[above]{$\simeq$} (-26) (26-) edge node[above]{$\epsilon$} (-27);
		\draw[->,font=\scriptsize] (14-) arc (-90:180:10pt); \draw[<-,font=\scriptsize] (-24) arc (-270:5:10pt);
	\end{tikzmath}
	where the first loop is
	\[
	p_\ast\pi_{1!}\pi_2^!p^\ast\stackrel c\simeq p_\ast f_! \pi_{1!}(f\times\id)^!\pi_2^!p^\ast\stackrel c\simeq 
	p_\ast \pi_{1!}(f\times\id)_!(f\times\id)^!\pi_2^!p^\ast\xrightarrow{\epsilon}p_\ast\pi_{1!}\pi_2^!p^\ast
	\]
	and the second loop is \[p^!\simeq \1_X\wedge p^!(\ph)\xrightarrow{\omega\wedge\id}\1_X\wedge p^!(\ph)\simeq p^!.\]
\end{proposition}

\begin{proof}
	By the base change theorem, the exchange transformations
	\[
	\Ex^!_!\colon \pi_{1!}\pi_2^!\to p^!p_!\quad\text{and}\quad \Ex_*^* \colon p^*p_*\to \pi_{2*}\pi_1^*
	\]
	are invertible. Lemma~\ref{lem:bc} shows that, under these isomorphisms, the first row of the given composition is the coevaluation~\eqref{eqn:coev}, and the second row is the evaluation~\eqref{eqn:ev}. Lemma~\ref{lem:symmetry} shows that the vertical arrow is inverse to the symmetry $p_*p^!\1_B\wedge p_! p^\ast\1_B\simeq p_! p^\ast\1_B\wedge p_*p^!\1_B$. It remains to prove that, under the isomorphism $\Sigma^\infty_+X\simeq p_*p^!\1_B$, the first loop corresponds to $\Sigma^\infty_+f\wedge\id$, and the second loop corresponds to $\id\wedge p_\sharp\omega$.
	
	Recall from~\S\ref{sec:review} that $\Sigma^\infty_+f$ is the following composition evaluated at $\1_B$:
	\[p_*p^!\stackrel c\simeq p_*f_!f^!p^!\stackrel\epsilon\to p_*p^!.\]
	Under the projection isomorphism~\eqref{eqn:duality2}, $\Sigma^\infty_+f\wedge\id$ is therefore the composition
	\begin{equation}\label{eqn:w78e}
		p_\ast p^!p_! p^\ast \simeq p_\ast f_! f^!p^!p_!p^\ast\xrightarrow{\epsilon} p_\ast p^!p_!p^\ast.
	\end{equation}
	Applying Lemma~\ref{lem:prism} to the pair of cartesian squares
	\begin{tikzmath}
		\diagram{X\times_BX & X\times_BX & X \\ X & X & B\rlap, \\};
		\arrows (11-) edge node[above]{$f\times \id$} (-12) (11) edge node[left]{$\pi_1$} (21) (21-) edge node[below]{$f$} (-22) (12-) edge node[above]{$\pi_2$} (-13) (12) edge node[left]{$\pi_1$} (22) (22-) edge node[below]{$p$} (-23) (13) edge node[right]{$p$} (23);
	\end{tikzmath}
	we deduce that~\eqref{eqn:w78e} becomes the first loop under the exchange isomorphisms $\Ex_!^!$.
	
	Denote also by $\omega$ the image of $\omega$ under the obvious map $\End(\1_X)\to \End(\id_{\SH(X)})$, so that the second loop is the natural transformation $p_!\pi_{2*}\pi_1^*\omega p^!$. By the compatibility of Thom transformations with the monoidal structure, the transformation $\omega$ commutes with any Thom transformation. The square
	\begin{tikzmath}
		\diagram{p_\sharp p^* & p_\sharp\Sigma^{-\Omega_p}\Sigma^{\Omega_p}p^* & p_*p^! \\ p_\sharp p^* & p_\sharp\Sigma^{-\Omega_p}\Sigma^{\Omega_p}p^* & p_*p^! \\};
		\arrows (11-) edge node[above]{$\simeq$} (-12) (12-) edge node[above]{$\simeq$} (-13) (21-) edge node[above]{$\simeq$} (-22) (22-) edge node[above]{$\simeq$} (-23)
		(11) edge node[left]{$p_\sharp\omega p^*$} (21) (13) edge node[right]{$p_*\omega p^!$} (23);
	\end{tikzmath}
	is therefore commutative. Under the natural isomorphism~\eqref{eqn:duality1}, $\id\wedge p_\sharp\omega$ then becomes $p_!p^*p_*\omega p^!$, which is the given loop modulo the exchange isomorphism $\Ex_*^*$.
\end{proof}

In the rest of this section we will give a more explicit description of this duality in a special case which will be used in \S\ref{sec:fields}. In what follows we often omit the stabilization functor $\Sigma^\infty$ from the notation and implictly view pointed presheaves on $\Sm_B$ as objects of $\SH(B)$ (we do not mean to say that the maps we consider are defined unstably, although this will sometimes be the case).

\begin{samepage}
\begin{definition}\label{def:euclidean}
	A \emph{Euclidean embedding} of $X$ is a triple $(s,V,\beta)$ where:
	\begin{itemize}
		\item $s$ is a closed immersion $X\into E$ in $\Sm_B$;
		\item $V$ is a vector bundle on $B$; \item $\beta$ is a path from $s^*(\bb V(\Omega_{E/B}))$ to $p^*(V)$ in the $K$-theory groupoid $K(X)$.
		\end{itemize}
\end{definition}
\end{samepage}

The proof of \cite[Lemma 2.8]{VV} shows that $X$ admits a Euclidean embedding if it is a closed subscheme of a projective bundle over $B$. Note also that, if $X$ admits a Euclidean embedding $(s,V,\beta)$, then it has one where the closed immersion is the zero section of a vector bundle, namely $X\into \bb V(\scr N_s)$. In addition to being smooth and proper, we now assume that $X$ admits a Euclidean embedding $(s,V,\beta)$, which we fix once and for all. The path $\beta$ in $K(X)$ determines an isomorphism
\begin{equation*}\label{eqn:Kpath}
	\Sigma^{s^\ast(\Omega_{E})}\simeq\Sigma^{p^*(V)}
\end{equation*}
of self-equivalences of $\SH(X)$.
The short exact sequence of locally free sheaves on $X$
\[0\to\scr N_s\to s^\ast(\Omega_{E})\xrightarrow{ds}\Omega_{X}\to 0\]
then induces an isomorphism
\[\Sigma^{-\Omega_{X}}\simeq\Sigma^{-p^*(V)} \Sigma^{\scr N_s},\]
whence
\begin{equation*}\label{eqn:stablenormal}
	p_!\simeq p_\sharp\Sigma^{-\Omega_{X}}\simeq p_\sharp \Sigma^{-p^*(V)} \Sigma^{\scr N_s}\simeq \Sigma^{-V}p_\sharp \Sigma^{\scr N_s}.
\end{equation*}
Finally, by the purity isomorphism, we obtain
\begin{equation}\label{eqn:dual}
	p_! \1_X\simeq \Sigma^{-V}\frac E{E\minus X}.
\end{equation}
It is worth emphasizing the the isomorphism~\eqref{eqn:dual} depends not only on $s$ and $V$ but also on $\beta$.

Under the isomorphism~\eqref{eqn:dual}, the coevaluation map~\eqref{eqn:coev} is the $V$-desuspension of a composition
\[
S^{V}\longrightarrow \frac{E}{E\minus X}\longrightarrow X_+\wedge \frac{E}{E\minus X},
\]
and the evaluation map~\eqref{eqn:ev} is the $V$-desuspension of a composition
\[
\frac{E}{E\minus X}\wedge X_+\longrightarrow \Sigma^{V}X_+\longrightarrow S^{V}.
\]
We would like to describe these four maps more explicitly. 
Let $\hat p\colon E\to B$ be the structure map of $E$, and define $\hat\pi_1$ and $\hat\pi_2$ by the cartesian square
\begin{tikzmath}
	\diagram{E\times X & X \\ E & B\rlap. \\};
	\arrows (11-) edge node[above]{$\hat\pi_2$} (-12) (11) edge node[left]{$\hat\pi_1$} (21) (21-) edge node[below]{$\hat p$} (-22) (12) edge node[right]{$p$} (22);
\end{tikzmath}
Let also $t\colon X\into E\times X$ be the composition $(s\times\id)\circ\delta$. We will define an isomorphism of short exact sequences
\begin{tikzequation}\label{eqn:SES}
	\diagram{0 & \scr N_s & s^*(\Omega_E) & \Omega_X & 0 \\
	0 & \delta^*(\scr N_{s\times\id}) & \scr N_t & \scr N_\delta & 0\rlap. \\};
	\arrows (11-) edge (-12) (12-) edge (-13) (13-) edge (-14) (14-) edge (-15)
	(21-) edge (-22) (22-) edge (-23) (23-) edge (-24) (24-) edge (-25)
	(12) edge node[left]{$\simeq$} node[right]{$\xi$} (22) (13) edge node[left]{$\simeq$} node[right]{$\mu$} (23) (14) edge node[right]{$\nu_2$} node[left]{$\simeq$} (24);
\end{tikzequation}
The isomorphism $\xi$ is the composition
\[\scr N_s\simeq\delta^*\pi_1^*(\scr N_s)\simeq \delta^*(\scr N_{s\times\id}).\]
The isomorphism $\nu_2\colon\Omega_X\stackrel\sim\to\scr N_\delta$ is defined so that the composition
\[\Omega_X\xrightarrow{\nu_2} \scr N_\delta\into\delta^*(\Omega_{X\times X})\]
is $\delta^*(d\pi_1)-\delta^*(d\pi_2)$.
In other words, $\nu_2$ is the composition of the canonical isomorphisms
\[\Omega_p\simeq \delta^*\pi_1^*(\Omega_p)\simeq \delta^*(\Omega_{\pi_2})\simeq\scr N_\delta.\]
It is then clear that the composite equivalence
\[\id\simeq \delta^!\pi_2^!\simeq \delta^!\Sigma^{\Omega_{\pi_2}}\pi_2^*\simeq \Sigma^{\delta^*(\Omega_{\pi_2})}\delta^!\pi_2^*\simeq \Sigma^{\Omega_p}\delta^!\pi_2^*\stackrel\Pi\simeq\Sigma^{\Omega_p-\scr N_\delta}\delta^*\pi_2^*\simeq\Sigma^{\Omega_p-\scr N_\delta}\]
is induced by $\nu_2$.
Finally, the isomorphism $\mu$ is defined so that the composition
\[s^*(\Omega_E)\stackrel\mu\to \scr N_t\into t^*(\Omega_{E\times X})\]
is $t^*(d\hat\pi_1)-t^*(d(s\hat\pi_2))$. It is easy to check that the diagram~\eqref{eqn:SES} commutes.

\begin{proposition}\label{prop:desc}
	Let $(s\colon X\into E,V,\beta)$ be a Euclidean embedding giving rise to the isomorphism~\eqref{eqn:dual}.
	\begin{enumerate}
		\item Suppose that $s\colon X\into E$ is the zero section of a vector bundle $r\colon E\to X$. Then the map $p_*p^*\1_B\stackrel\eta\to p_*p^!p_!p^*\1_B\simeq p_*p^!\1_B\wedge p_!p^*\1_B$ is the $V$-desuspension of the composition
		\begin{equation*}
			\frac{E}{E\minus X}\xrightarrow{(r,\id)}\frac{X\times E}{X\times(E\minus X)}\simeq X_+\wedge \frac{E}{E\minus X}.
		\end{equation*}
		\item The map $p_!p^*\1_B\wedge p_*p^!\1_B\simeq p_!p^*p_*p^!\1_B\stackrel\epsilon\to p_!p^!\1_B$ is the $V$-desuspension of the composition
		\begin{equation*}
			\frac{E}{E\minus X}\wedge X_+\simeq\frac{E\times X}{(E\minus X)\times X} \to \frac{E\times X}{(E\times X)\minus \Delta_X}\stackrel\Pi\simeq \Th_X(\bb V(\scr N_t))\simeq \Sigma^{V}X_+,
		\end{equation*}
		where the last isomorphism is induced by $\mu\colon \scr N_t\simeq s^*(\Omega_E)$ and by $\beta$.
	\end{enumerate}
\end{proposition}

\begin{proof}
	(1) We must check that the two outer compositions in the following diagram coincide:
	\begin{tikzmath}
		\def\colsep{1.8em}
		\diagram{p_*\Sigma^Ep^!\1_B & p_*\pi_{1!}\delta_!\delta^!\pi_2^!\Sigma^Ep^!\1_B & p_*\Sigma^Ep^!\1_B & \\
		& p_*\pi_{1!}\pi_2^!\Sigma^Ep^!\1_B & p_*p^!p_!\Sigma^Ep^!\1_B & p_*p^!\1_B\wedge p_!\Sigma^Ep^!\1_B \\
		& p_\sharp\pi_{1\sharp}\pi_2^*\Sigma^Ep^*\1_B & p_\sharp p^*p_\sharp\Sigma^Ep^*\1_B & p_\sharp p^*\1_B\wedge p_\sharp\Sigma^Ep^*\1_B \\
		\Th_X(E) & \Th_{X\times X}(\pi_2^*E) & & X_+\wedge\Th_X(E)\rlap. \\};
		\arrows
		(11-) edge node[below]{$\simeq$} (-12) (-13) edge node[below]{$\simeq$} (12-)
		(11) edge node[left]{$\simeq$} (41) (12) edge node[left]{$\epsilon$} (22) (13) edge node[right]{$\eta$} (23)
		(41-) edge node[above]{$(r,\id)$} (-42)
		(22-) edge node[above]{$\Ex_!^!$} (-23) (23-) edge[<-] node[above]{$\Pr_!^*$} (-24)
		(32-) edge node[above]{$\Ex_\sharp^*$} (-33) (33-) edge node[above]{$\Pr_\sharp^*$} (-34)
		(42-) edge node[above]{$\simeq$} (-44)
		(22) edge node[left]{$\simeq$} (32) (23) edge node[right]{$\simeq$} (33) (24) edge node[right]{$\simeq$} (34)
		(32) edge node[left]{$\simeq$} (42) (34) edge node[right]{$\simeq$} (44);
		\draw[->,font=\scriptsize] (11) to[out=15,in=165] node[above]{$\id$} (13);
	\end{tikzmath}
	The three vertical isomorphisms in the second row are obtained by getting rid of shrieks and rearranging the resulting Thom transformations. Note that $(r,\id)\colon \Th_X(E)\to \Th_{X\times X}(\pi_2^*E)$ is the map induced by the diagonal $\delta\colon X\to X\times_BX$ and the canonical isomorphism $\delta^*\pi_2^*E\simeq E$. We saw in~\S\ref{sec:review} that the left-hand rectangle is commutative. The commutativity of the top square is Lemma~\ref{lem:bc} (2). Finally, one verifies easily that the lower rectangle is the stabilization of a commutative rectangle of presheaves of pointed sets on $\Sm_B$. Thus, the whole diagram is commutative.
	
	(2) We first express the given composition in terms of the six operations. We have
	\[\frac{E\times X}{(E\minus X)\times X}\simeq p_\sharp \hat\pi_{2\sharp}(s\times\id)_*\pi_1^*p^*\1_B
	\quad\text{and}\quad
	\frac{E\times X}{(E\times X)\minus\Delta_X}\simeq p_\sharp\hat\pi_{2\sharp}t_*p^*\1_B,\]
	and the map
	\[\frac{E\times X}{(E\minus X)\times X}\to \frac{E\times X}{(E\times X)\minus\Delta_X}\]
	collapsing the complement of the diagonal is given by
	\begin{equation}\label{eqn:ev1}
		p_\sharp \hat\pi_{2\sharp}(s\times\id)_*\pi_1^*p^*\1_B\stackrel\eta\to p_\sharp \hat\pi_{2\sharp}(s\times\id)_*\delta_*\delta^*\pi_1^*p^*\1_B\simeq p_\sharp\hat\pi_{2\sharp}t_*p^*\1_B,
	\end{equation}
	as one can see at the level of pointed presheaves on $\Sm_B$.
	Consider the following diagram:
	\begin{tikzequation}\label{eqn:ev2}
		\def\colsep{1.8em}
		\diagram{p_\sharp\Sigma^{\scr N_s}p^*\1_B\wedge p_\sharp p^*\1_B & p_\sharp p^* p_\sharp \Sigma^{\scr N_s}p^*\1_B & p_\sharp \pi_{2\sharp}\pi_1^* \Sigma^{\scr N_s}p^*\1_B & p_\sharp \pi_{2\sharp}\Sigma^{\pi_1^*\scr N_s}\pi_1^*p^*\1_B & \\
		\hat p_\sharp s_* p^*\1_B\wedge p_\sharp p^*\1_B & p_\sharp p^*\hat p_\sharp s_* p^*\1_B & p_\sharp \hat\pi_{2\sharp}\hat \pi_1^* s_* p^*\1_B & p_\sharp\hat\pi_{2\sharp}(s\times\id)_*\pi_1^*p^*\1_B \\
		\frac E{E\minus X}\wedge X_+ & & & \frac{E\times X}{(E\minus X)\times X}\rlap. \\
		};
		\arrows (11-) edge[<-] node[above]{$\Pr_\sharp^*$} (-12) (12-) edge[<-] node[above]{$\Ex_\sharp^*$} (-13) (13-) edge node[above]{$\simeq$} (-14)
		(21-) edge[<-] node[above]{$\Pr_\sharp^*$} (-22) (22-) edge[<-] node[above]{$\Ex_\sharp^*$} (-23) (23-) edge node[above]{$\Ex^*_*$} (-24) (11) edge node[left]{$\Pi$} (21) (12) edge node[left]{$\Pi$} (22) (14) edge node[left]{$\Pi$} (24)
		(31-) edge node[above]{$\simeq$} (-34) (21) edge node[left]{$\simeq$} (31) (24) edge node[right]{$\simeq$} (34);
	\end{tikzequation}
	The lower rectangle is seen to be commutative at the level of pointed presheaves, and the top left square commutes by naturality of $\Pr_\sharp^*$. The upper right rectangle becomes an instance of the compatibility of $\Ex_*^*$ with compositions of cartesian squares after replacing lower sharps by lower stars.
	It remains to prove that the composition
	\[
	p_\sharp p^*p_*\Sigma^{s^*(\Omega_E)}p^*\simeq p_\sharp p^* p_\sharp \Sigma^{\scr N_s}p^*\xrightarrow{\eqref{eqn:ev2}} p_\sharp\hat\pi_{2\sharp}(s\times\id)_*\pi_1^*p^* \xrightarrow{\eqref{eqn:ev1}} p_\sharp\hat\pi_{2\sharp}t_*p^*\stackrel\Pi\simeq p_\sharp\Sigma^{\scr N_t}p^*\stackrel\mu\simeq p_\sharp\Sigma^{s^*(\Omega_E)}p^*
	\]
	is equal to the counit $\epsilon$ (when evaluated at $\1_B$). One finds these two maps as the boundary of the following diagram, after applying $p_\sharp(\ph)p^*$:
	\begin{tikzmath}
		\diagram{p^*p_*\Sigma^{s^*(\Omega_E)} & \pi_{2*}\pi_1^* \Sigma^{s^*(\Omega_E)} & \pi_{2*}\delta_*\Sigma^{s^*(\Omega_E)} & \Sigma^{s^*(\Omega_E)} \\
		p^*p_\sharp\Sigma^{\scr N_s} & \pi_{2\sharp}\pi_1^*\Sigma^{\scr N_s} & \pi_{2\sharp}\delta_*\Sigma^{\scr N_s} &  \\
		& \pi_{2\sharp}\Sigma^{\scr N_{s\times\id}}\pi_1^* & \pi_{2\sharp}\Sigma^{\scr N_{s\times\id}}\delta_* & \Sigma^{\scr N_t} \\
		& \hat\pi_{2\sharp}(s\times\id)_*\pi_1^* & \hat\pi_{2\sharp}(s\times\id)_*\delta_* & \hat\pi_{2\sharp}t_*\rlap. \\};
		\arrows (11-) edge node[below]{$\Ex_*^*$} (-12) (12-) edge node[below]{$\eta$} (-13) (13-) edge node[below]{$\simeq$} (-14)
		(11) edge node[left]{$\simeq$} (21) (12) edge node[left]{$\simeq$} (22) (13) edge node[left]{$\simeq$} (23)
		(21-) edge[<-] node[below]{$\Ex_\sharp^*$} (-22) (22-) edge node[below]{$\eta$} (-23)
		(22) edge node[left]{$\simeq$} (32) (23) edge node[left]{$\simeq$} (33) (32-) edge node[below]{$\eta$} (-33)
		(32) edge[<-] node[left]{$\Pi$} (42) (33) edge[<-] node[left]{$\Pi$} (43) (42-) edge node[below]{$\eta$} (-43)
		(43-) edge node[below]{$\simeq$} (-44) (44) edge node[right]{$\Pi$} (34) (34) edge[<-] node[right]{$\mu$} (14)
		;
		\draw[->,font=\scriptsize] (11) to[out=15,in=165] node[above]{$\epsilon$} (14);
	\end{tikzmath}
	We claim that this diagram commutes. The topmost face commutes by Lemma~\ref{lem:bc} (1), and the commutativity of the four small squares is clear. The large rectangle may be decomposed as follows:
	\begin{tikzmath}
		\diagram{ & \pi_{2*}\delta_* \Sigma^{s^*(\Omega_E)} & \Sigma^{s^*(\Omega_E)} \\
		\pi_{2\sharp}\delta_*\Sigma^{\scr N_s} & \pi_{2*}\delta_*\Sigma^{\Omega_X}\Sigma^{\scr N_s} & \Sigma^{\Omega_X}\Sigma^{\scr N_s} \\
		\pi_{2\sharp}\Sigma^{\scr N_{s\times\id}}\delta_* & \pi_{2\sharp}\delta_* \Sigma^{\delta^*(\scr N_{s\times\id})} & \Sigma^{\scr N_\delta}\Sigma^{\delta^*(\scr N_{s\times\id})} \\
		\hat\pi_{2\sharp}(s\times\id)_*\delta_* & \hat\pi_{2\sharp}t_* & \Sigma^{\scr N_t}\rlap. \\};
		\arrows 
		(12-) edge node[above]{$\simeq$} (-13)
		(12) edge node[above left]{$\simeq$} (21)
		(12) edge node[right]{$\simeq$} (22) (13) edge node[left]{$\simeq$} (23)
		(21-) edge node[below]{$\simeq$} (-22) (22-) edge node[below]{$\simeq$} (-23)
		(21) edge node[left]{$\simeq$} (31) (23) edge node[left]{$\nu_2$} node[right]{$\xi$} (33)
		(31-) edge node[above]{$\simeq$} (-32) (32-) edge node[above]{$\Pi$} (-33)
		(31) edge[<-] node[left]{$\Pi$} (41) (33) edge[<-] node[left]{$\simeq$} (43)
		(41-) edge node[above]{$\simeq$} (-42) (42-) edge node[above]{$\Pi$} (-43)
		;
		\draw[<-,font=\scriptsize] ([xshift=3pt] 43.north east) to[out=45,in=315] node[right]{$\mu$} ([xshift=-4pt] 13.south east);
	\end{tikzmath}
	The rightmost face commutes by~\eqref{eqn:SES} and the middle rectangle commutes by the definitions of $\xi$ and $\nu_2$. Finally, the bottom rectangle commutes by the compatibility of the purity isomorphisms with the composition of the closed immersions $\delta$ and $s\times\id$ \cite[\S1.6.4]{Ayoub}.
	 \end{proof}

The counit $\epsilon\colon p_!p^!\1_B\to \1_B$ is of course the map $\Sigma^\infty_+ p\colon\Sigma^\infty_+ X\to \Sigma^\infty_+ B=\1_B$.
The unit $\eta\colon \1_B\to p_*p^*\1_B$ is more difficult to describe explicitly, and we do not know how to do it in any kind of generality.\footnote{If $X$ is a closed subscheme of a projective bundle over $B$, it is possible that the unstable map $S^{V}\to E/(E\minus X)$ constructed in \cite[Theorem 2.11 (2)]{VV} (for a specific Euclidean embedding $(s,V,\beta)$) stabilizes to $\eta$, but we did not check it.} However, we can at least give a useful characterization of $\eta$:

\begin{proposition}\label{prop:eta}
	Let $(s\colon X\into E,V,\beta)$ be a Euclidean embedding and let
	\[\zeta\colon S^{V}\to \frac{E}{E\minus X}\]
	be a map in $\SH(B)$. The following conditions are equivalent:
	\begin{enumerate}
		\item Under the isomorphism~\eqref{eqn:dual}, $\zeta$ is the $V$-suspension of the unit $\eta\colon\1_B\to p_\ast p^\ast \1_B$.
		\item The composition
\[\Sigma^{V} X_+\xrightarrow{\zeta\wedge\id}\frac{E}{E\minus X}\wedge X_+\to \Sigma^{V} X_+\stackrel p\to S^{V},\]
where the second map is that given in Proposition~\ref{prop:desc} (2), is equal to $\Sigma^{V}\Sigma^\infty_+p$.
	\end{enumerate}
\end{proposition}

\begin{proof}
	Since the unit $\eta\colon\1_B\to p_\ast p^\ast \1_B$ is dual to $\epsilon\colon p_!p^!\1_B\to \1_B$, it is determined by the equation $\ev\circ(\eta\wedge\id)=\epsilon$. The equivalence of (1) and (2) is now clear by Proposition~\ref{prop:desc} (2).
	\end{proof}

In \S\ref{sec:fields}, we will define a map $\zeta$ satisfying the condition of Proposition~\ref{prop:eta} (2) when $B$ is a field and $X$ is a finite separable extension of $B$. As a result, the duality will be completely explicit in this case.

\newpage
\section{Proof of the main theorem}
\label{sec:traces}

We prove Theorem~\ref{thm:lefschetz}. As a warm-up, assume that $X$ admits a Euclidean embedding $(s,V,\beta)$, chosen such that $s$ is the zero section of a vector bundle $r\colon E\to X$.
By Proposition~\ref{prop:desc}, the trace of $\Sigma^\infty_+f$ is then the $V$-desuspension of a composition
\begin{multline*}
S^{V}\to \frac{E}{E\minus X}\xrightarrow{(r,\id)}\frac{X\times E}{X\times (E\minus X)}\xrightarrow{f\times \id}
\frac{X\times E}{X\times (E\minus X)}\\
\stackrel{\tau}{\simeq}\frac{E\times X}{(E\minus X)\times X}\to\frac{E\times X}{(E\times X)-\Delta_X}\simeq \Sigma^{V}X_+\to S^{V},
\end{multline*}
where $\tau$ is the monoidal symmetry. Ignoring the first and last arrows, it is clear that the remaining composition factors through $E/(E\minus X^f)$, and hence that $\tr(\Sigma^\infty_+ f)=0$ if $X^f$ is empty (compare this argument with the proof of the Lefschetz–Hopf theorem in \cite{DP}). It is possible to prove the more precise statement of Theorem~\ref{thm:lefschetz} in this explicit setting, but, to treat the general case where $X$ is proper over an arbitrary base, we will now switch to the formalism of six operations.

Throughout this section we use the following notation:
\[
\begin{tikzpicture}
	\diagram{X^f & X \\ & B\rlap, \\};
	\arrows (11-) edge[c->] node[above]{$i$} (-12) (12) edge node[right]{$p$} (22) (11) edge node[below left]{$q$} (22);
\end{tikzpicture}
\qquad
\begin{tikzpicture}
	\diagram{X\times_BX & X \\ X & B\rlap, \\};
	\arrows (11-) edge node[above]{$\pi_2$} (-12) (11) edge node[left]{$\pi_1$} (21) (21-) edge node[below]{$p$} (-22) (12) edge node[right]{$p$} (22);
\end{tikzpicture}
\qquad
\begin{tikzpicture}
	\diagram{X^f\times_BX^f & X^f \\ X^f & B\rlap. \\};
	\arrows (11-) edge node[above]{$\rho_2$} (-12) (11) edge node[left]{$\rho_1$} (21) (21-) edge node[below]{$q$} (-22) (12) edge node[right]{$q$} (22);
\end{tikzpicture}
\]
Let also $\delta\colon X\into X\times_BX$ and $\theta\colon X^f\into X^f\times_BX^f$ be the diagonal maps, and let $\gamma=(f\times\id)\circ\delta$ be the graph of $f$.
For the moment we do not assume that $f$ has regular fixed points.

Recall from Proposition~\ref{prop:trace} that $\tr(\Sigma^\infty_+f)$ can be expressed as a certain composition
\[\id\stackrel\eta\to p_*p^*\to p_!p^!\stackrel\epsilon\to \id\]
evaluated at $\1_B$, where the transformation $p_*p^*\to p_!p^!$ is the solid boundary of the following commutative diagram:
\begin{tikzequation}\label{eqn:magic}
	\diagram{
	p_\ast p^\ast & p_\ast\pi_{1!} \delta_!\delta^! \pi_2^!p^\ast & p_\ast \pi_{1!}\pi_2^!p^\ast \\
	& p_\ast\pi_{1!}(f\times\id)_!\delta_!\delta^!(f\times\id)^!\pi_2^!p^\ast & p_\ast\pi_{1!}(f\times\id)_!(f\times\id)^!\pi_2^!p^\ast \\
	& p_\ast\pi_{1!}\gamma_!\gamma^!\pi_2^!p^\ast & p_\ast\pi_{1!}\pi_2^!p^\ast \\
	& p_!\pi_{2\ast}\gamma_!\gamma^!\pi_1^\ast p^! & p_!\pi_{2\ast}\pi_1^\ast p^! \\
	& p_!\pi_{2\ast}\gamma_!\gamma^!\delta_\ast\delta^\ast\pi_1^\ast p^! & p_!\pi_{2\ast}\delta_\ast\delta^\ast\pi_1^\ast p^! \\
	p_! i_\ast i^! p^! & p_!\gamma^!\delta_\ast p^! & p_! p^!\rlap. \\
	};
	\arrows (11) edge[dashed] (61)
	(11-) edge node[above]{$\simeq$} (-12)
	(12-) edge node[above]{$\epsilon$} (-13)
	(22-) edge node[above]{$\epsilon$} (-23)
	(32-) edge node[above]{$\epsilon$} (-33)
	(42-) edge node[above]{$\epsilon$} (-43)
	(52-) edge node[above]{$\epsilon$} (-53)
	(62-) edge[dashed] (-63)
	(61-) edge node[below]{$\simeq$} node[above]{$\Ex^!_\ast$} (-62)
	(12) edge node[left]{$\simeq$} (22)
	(13) edge node[right]{$\simeq$} (23)
	(22) edge node[left]{$\simeq$} (32)
	(23) edge node[right]{$\epsilon$} (33)
	(32) edge node[left]{$\simeq$} (42)
	(33) edge node[left]{$\simeq$} node[right]{$\tau$} (43)
	(42) edge node[left]{$\eta$} (52)
	(43) edge node[right]{$\eta$} (53)
	(52) edge node[left]{$\simeq$} (62)
	(53) edge node[right]{$\simeq$} (63)
	;
\end{tikzequation}
The map at the bottom left is the exchange transformation $\Ex^!_\ast\colon i_\ast i^!\to\gamma^!\delta_\ast$ associated with the cartesian square
\begin{tikzequation}\label{eqn:fixedpointsquare}
	\diagram{X^f & X \\ X & X\times_BX\rlap, \\};
	\arrows (11-) edge node[above]{$i$} (-12) (11) edge node[left]{$i$} (21) (21-) edge node[below]{$\gamma$} (-22) (12) edge node[right]{$\delta$} (22);
\end{tikzequation}
and it is an isomorphism by the base change theorem. The dashed arrows in~\eqref{eqn:magic} can then be defined so as to make the diagram commute. Lemma~\ref{lem:bottomrow} shows that the bottom row in~\eqref{eqn:magic} is the counit $\epsilon\colon p_! i_\ast i^! p^!\to p_!p^!$. Note that this diagram already proves Corollary~\ref{cor:lefschetz}: if $X^f$ is empty, then $i_\ast i^!$ is the zero functor and hence $\tr(\Sigma^\infty_+f)=0$.

The dashed arrow $p_\ast p^\ast\to p_!i_\ast i^!p^!$ in~\eqref{eqn:magic} is the composition of the right column and the bottom row in the following diagram:
\begin{tikzequation}\label{eqn:magic2}
	\diagram{
	p_\ast i_\ast i^\ast p^\ast & p_\ast p^\ast \\
	p_\ast\pi_{1!} \delta_!i_\ast i^\ast\delta^! \pi_2^!p^\ast & p_\ast\pi_{1!} \delta_!\delta^! \pi_2^!p^\ast \\
	p_\ast\pi_{1!}\gamma_!i_\ast i^\ast\gamma^!\pi_2^!p^\ast & p_\ast\pi_{1!}\gamma_!\gamma^!\pi_2^!p^\ast \\
	p_!\pi_{2\ast}\gamma_!i_\ast i^\ast\gamma^!\pi_1^\ast p^! & p_!\pi_{2\ast}\gamma_!\gamma^!\pi_1^\ast p^! \\
	p_!\pi_{2\ast}\gamma_!i_\ast i^!\delta^\ast\pi_1^\ast p^! & p_!\pi_{2\ast}\gamma_!\gamma^!\delta_\ast\delta^\ast\pi_1^\ast p^! \\
	p_! i_\ast i^! p^! & p_!\gamma^!\delta_\ast p^!\rlap. \\
	};
	\arrows
	(-12) edge node[above]{$\eta$} (11-)
	(-22) edge node[above]{$\eta$} (21-)
	(-32) edge node[above]{$\eta$} (31-)
	(-42) edge node[above]{$\eta$} (41-)
	(51-) edge node[below]{$\simeq$} node[above]{$\Ex^!_\ast$} (-52)
	(61-) edge node[below]{$\simeq$} node[above]{$\Ex^!_\ast$} (-62)
	(11) edge node[left]{$\simeq$} (21)
	(21) edge node[left]{$\simeq$} (31)
	(31) edge node[left]{$\simeq$} (41)
	(41) edge node[left]{$\Ex^{\ast !}$} (51)
	(51) edge node[left]{$\simeq$} (61)
	(12) edge node[right]{$\simeq$} (22)
	(22) edge node[right]{$\simeq$} (32)
	(32) edge node[right]{$\simeq$} (42)
	(42) edge node[right]{$\eta$} (52)
	(52) edge node[right]{$\simeq$} (62)
	;
\end{tikzequation}
The commutativity of each square in this diagram is clear, except that of the fourth square which follows from the definition of the exchange transformation $\Ex_*^!$ in terms of $\Ex^{*!}$ (see \cite[\S1.2.4]{Ayoub}).

The left column of~\eqref{eqn:magic2} is a natural transformation $q_\ast q^\ast \to q_!q^!$ which, by~\eqref{eqn:magic} and Lemma~\ref{lem:bottomrow}, makes the following diagram commute:
\begin{tikzequation}\label{eqn:diamond}
	\def\rowsep{1em}
	\diagram{ & p_\ast p^\ast & p_!p^! & \\ \id & & & \id\rlap. \\ & q_\ast q^\ast & q_! q^! & \\};
	\arrows (12-) edge (-13) (32-) edge node[above]{\eqref{eqn:magic2}} (-33) (12) edge node[left]{$\eta$} (32) (33) edge node[right]{$\epsilon$} (13)
	(21) edge node[above]{$\eta$} (12) edge node[below]{$\eta$} (32) (13) edge node[above]{$\epsilon$} (24) (33) edge node[below]{$\epsilon$} (24);
\end{tikzequation}

Assume now that $f$ has regular fixed points, \ie, that $q$ is smooth and that $\id-i^\ast(df)$ restricts to an automorphism $\phi\colon \scr N_i\stackrel\sim\to\scr N_i$. By Proposition~\ref{prop:trace}, $\int_{X^f}\J\phi\d\chi$ is a certain composition
\begin{equation}\label{eqn:tracethom}
	\id\stackrel\eta\to q_* q^*\to q_!q^!\stackrel\epsilon\to \id
\end{equation}
evaluated at $\1_B$. In view of~\eqref{eqn:diamond}, to conclude the proof of Theorem~\ref{thm:lefschetz}, it will suffice to show that the segment $q_\ast q^\ast \to q_!q^!$ in~\eqref{eqn:tracethom} is equal to the left column of~\eqref{eqn:magic2}. This segment (as given by Proposition~\ref{prop:trace}) is the composition of the top row and the right vertical arrows in the following diagram:
\begin{tikzequation}\label{eqn:twistedEuler}
	\diagram{
	q_\ast q^\ast & q_\ast\rho_{1!} \theta_!\theta^! \rho_2^!q^\ast & q_\ast \rho_{1!}\rho_2^!q^\ast \\
	& q_!\rho_{2\ast}\theta_!\theta^!\rho_1^\ast q^! & q_!\rho_{2\ast}\rho_1^\ast q^! \\
	& q_!\rho_{2\ast}\theta_!\theta^!\rho_1^\ast q^! & q_!\rho_{2\ast}\rho_1^\ast q^! \\
	& q_!\rho_{2\ast}\theta_!\theta^!\theta_\ast\theta^\ast\rho_1^\ast q^! & q_!\rho_{2\ast}\theta_\ast\theta^\ast\rho_1^\ast q^! \\
	& q_!\rho_{2\ast}\theta_!\theta^\ast\rho_1^\ast q^! & q_! q^!\rlap, \\
	};
	\arrows (11-) edge node[above]{$\simeq$} (-12)
	(12) edge node[right]{$\simeq$} (22)
	(22) edge node[right]{$\simeq$} node[left]{$\J\phi$} (32)
	(32) edge node[left]{$\eta$} (42)
	(42) edge[<-] node[right]{$\simeq$} node[left]{$\eta$} (52)
	(13) edge node[left]{$\simeq$} node[right]{$\tau$} (23)
	(23) edge node[left]{$\simeq$} node[right]{$\J\phi$} (33)
	(33) edge node[right]{$\eta$} (43)
	(43) edge node[left]{$\simeq$} (53)
	(12-) edge node[above]{$\epsilon$} (-13)
	(22-) edge node[above]{$\epsilon$} (-23)
	(32-) edge node[above]{$\epsilon$} (-33)
	(42-) edge node[above]{$\epsilon$} node[below]{$\simeq$} (-43)
	(52-) edge node[above]{$\simeq$} (-53) 
	;
	\draw[->,font=\scriptsize] (-32) to[out=225,in=135] node[left]{$\sigma$} (52.west);
\end{tikzequation}
where $\J\phi$ acts after $q^!$ in both columns. Each square commutes by the naturality of the given transformations, except the last square which commutes by a triangle identity for the adjunction $\theta_*\simeq\theta_!\dashv \theta^!$. The triangle at the bottom left commutes by Lemma~\ref{lem:sigma}.

\begin{remark}\label{rmk:vanishing}
	In the diagram~\eqref{eqn:twistedEuler}, one can replace $\J\phi$ by any endomorphism of $\1_{X^f}$ and $X^f$ itself by any smooth proper $B$-scheme. Theorem~\ref{thm:vanishing} follows from the observation that the natural transformation $\sigma\colon \theta^!\rho_1^\ast\to\theta^\ast\rho_1^\ast$ is zero if $[\Omega_{X^f}]=[\scr O_{X^f}]+[\scr E]$ in $K_0(X^f)$. Indeed, by Proposition~\ref{prop:purity}, this transformation can be identified with the transformation $\Sigma^{-\scr N_\theta}\to\Sigma^0$ induced by the epimorphism $\scr N_\theta\to 0$. Since $\scr N_\theta$ is isomorphic to $\Omega_{X^f}$, this transformation factors through the transformation $\Sigma^{-\A^1}\to\Sigma^{0}$ induced by the zero section of the trivial line bundle, which is clearly zero (see for example \cite[Lemme 1.6.1]{Ayoub}).
\end{remark}

We now compare the left column of~\eqref{eqn:magic2} with the left column of~\eqref{eqn:twistedEuler}. Both columns are of the form $q_*(q^*\to q^!)$, where the respective maps $q^*\to q^!$ are the left and right columns of the following diagram:
\begin{tikzmath}
	\diagram{
	& q^* & \\
	i^\ast\delta^!\pi_2^!p^\ast & & \theta^! \rho_2^!q^\ast \\
	i^\ast\gamma^!\pi_2^!p^\ast & & \theta^! \rho_2^!q^\ast \\
	i^\ast\gamma^!\pi_1^\ast p^! & & \theta^!\rho_1^\ast q^! \\
	& & \theta^!\rho_1^\ast q^! \\
	i^!\delta^\ast\pi_1^\ast p^! & & \theta^\ast\rho_1^\ast q^! \\
	& q^!\rlap. & \\
	};
	\arrows
	(21) edge node[right]{$\simeq$} node[left]{$c$} (31)
	(23) edge[-,vshift=1pt] (33) edge[-,vshift=-1pt] (33)
	(31) edge[<-] node[right]{$\simeq$} node[left]{$\Ex^{*!}$} (41)
	(41) edge node[left]{$\Ex^{\ast !}$} (61)
	(21-) edge node[below]{$\simeq$} node[above]{$\Ex^{*!}$} (-23)
	(31-) edge node[below]{$\simeq$} node[above]{$\Ex^{*!}$} (-33)
	(41-) edge node[below]{$\simeq$} node[above]{$\alpha$} (-43)
	(61-) edge[<-] node[above]{$\simeq$} node[below]{$\Ex^{*!}$} (-63)
	(33) edge[<-] node[left]{$\simeq$} node[right]{$\Ex^{*!}$} (43)
	(43) edge node[right]{$\J\phi$} node[left]{$\simeq$} (53)
	(53) edge node[right]{$\sigma$} (63)
	(12) edge node[above left]{$\simeq$} (21) edge node[above right]{$\simeq$} (23)
	(61) edge node[below left]{$\simeq$} (72) (63) edge node[below right]{$\simeq$} (72)
	;
\end{tikzmath}
Here the isomorphism $\alpha$ is defined by the commutativity of the second square. The commutativity of the first square is clear. Theorem~\ref{thm:lefschetz} is thus reduced to the commutativity of the pentagon
\begin{tikzmath}
	\diagram{i^\ast\gamma^! \pi_1^\ast p^! & \theta^!\rho_1^\ast q^! \\
	i^!\delta^\ast \pi_1^\ast p^! & \theta^\ast\rho_1^\ast q^!\rlap. \\};
	\arrows (11-) edge node[above]{$\simeq$} node[below]{$\alpha$} (-12) (11) edge node[left]{$\Ex^{\ast !}$} (21) (21-) edge node[above]{$\simeq$} (-22) (12) edge node[right]{$\sigma$} (22);
	\draw[<-,font=\scriptsize] (12-) arc (-90:185:10pt) node[right=28pt,above=0pt]{$\J\phi$};
\end{tikzmath}
This is the heart of the proof. By transforming the stars into shrieks, this pentagon becomes
\begin{tikzequation}\label{eqn:key}
	\diagram{\Sigma^{\scr N_i-i^\ast(\Omega_p)}q^! & \Sigma^{-\Omega_q}q^! \\
	\Sigma^{i^\ast(\scr N_\delta)-i^\ast(\Omega_p)}q^! & \Sigma^{\scr N_\theta-\Omega_q}q^!\rlap, \\};
	\arrows (11-) edge node[above]{$\simeq$} (-12) (11) edge (21) (21-) edge node[above]{$\simeq$} (-22) (12) edge (22);
	\draw[<-,font=\scriptsize] (12-) arc (-90:185:10pt) node[right=28pt,above=0pt]{$\J\phi$};
\end{tikzequation}
and we now identify the four unlabeled arrows. By definition of $\alpha$, the top map in~\eqref{eqn:key} is induced by the short exact sequence
\[0\to\scr N_i\to i^\ast(\Omega_p)\xrightarrow{di} \Omega_q\to 0.\]
 Denote by $\nu_1\colon\Omega_X\stackrel\sim\to\scr N_\delta$ the isomorphism for which the composition
\[\Omega_X\xrightarrow{\nu_1}\scr N_\delta\into \delta^\ast(\Omega_{X\times_BX})\]
is $\delta^\ast(d\pi_2)-\delta^\ast(d\pi_1)$. The composite isomorphism
\[\id\simeq \delta^\ast\pi_1^\ast\simeq\Sigma^{\scr N_\delta-\Omega_p}\delta^!\pi_1^!\simeq \Sigma^{\scr N_\delta-\Omega_p}\]
is then induced by $\nu_1$, and similarly for the isomorphism $\id\simeq \Sigma^{\scr N_\theta-\Omega_{q}}$ (for more details, see the discussion of the isomorphism $\nu_2$ before Proposition~\ref{prop:desc}). Under these trivializations, the bottom map in~\eqref{eqn:key} is just the identity $q^!\to q^!$.
The vertical maps in~\eqref{eqn:key} can be identified using Proposition~\ref{prop:purity}. Applying Proposition~\ref{prop:purity} to the cartesian square~\eqref{eqn:fixedpointsquare} shows that the left vertical arrow in~\eqref{eqn:key} is $\Sigma^{\psi}$ where $\psi\colon i^\ast(\scr N_\delta)\onto \scr N_i$ is the epimorphism induced by~\eqref{eqn:fixedpointsquare}. Explicitly, $\psi$ is determined by the following diagram of short exact sequences:
\begin{tikzmath}
	\diagram{i^\ast(\scr N_\delta) & i^\ast\delta^\ast(\Omega_{X\times_BX}) & i^\ast(\Omega_X) \\
	& i^\ast\gamma^\ast(\Omega_{X\times_BX}) & \\
	\scr N_i & i^\ast(\Omega_X) & \Omega_{X^f}\rlap. \\};
	\arrows (11-) edge[c->] (-12) (11) edge[dashed] node[left]{$\psi$} (31) (22) edge node[left]{$i^\ast(d\gamma)$} (32)
	(31-) edge[c->] (-32) (12-) edge[->>] node[above]{$i^\ast(d\delta)$} (-13) (32-) edge[->>] node[below]{$di$} (-33) (13) edge node[right]{$di$} (33) (12) edge node[left]{$\simeq$} (22);
\end{tikzmath}
Finally, applying Proposition~\ref{prop:purity} to the pullback of $\theta$ along itself shows that the right arrow in~\eqref{eqn:key} is $\Sigma^\zeta$ where $\zeta$ is the epimorphism $\scr N_\theta\to 0$. The commutativity of~\eqref{eqn:key} is thereby reduced to the commutativity of the following diagram:
\begin{tikzmath}
	\diagram{\Sigma^{\scr N_i -i^\ast(\Omega_X)} & \Sigma^{\scr N_i - i^\ast(\Omega_X)} & \Sigma^{-\Omega_{X^f}} \\
	\Sigma^{i^\ast(\scr N_\delta)-i^\ast(\Omega_X)} & \Sigma^{0} & \Sigma^{\scr N_\theta - \Omega_{X^f}}\rlap. \\};
	\arrows (11-) edge node[above]{$\phi$} node[below]{$\simeq$} (-12) (11) edge node[left]{$\psi$} (21) (12-) edge node[above]{$\simeq$} (-13)
	(21-) edge node[below]{$i^\ast(\nu_1)$} node[above]{$\simeq$} (-22) (22-) edge[<-] node[below]{$\nu_1$} node[above]{$\simeq$} (-23) (13) edge node[right]{$\zeta$} (23);
\end{tikzmath}
Applying $\Sigma^{i^*(\Omega_X)}$, it is equivalent to check that the following diagram commutes:
\begin{tikzmath}
	\diagram{\Sigma^{\scr N_i} & \Sigma^{\scr N_i} & \Sigma^{\scr N_\theta + \scr N_i} \\
	\Sigma^{i^\ast(\scr N_\delta)} & \Sigma^{i^\ast(\Omega_X)} & \Sigma^{\Omega_{X^f} +\scr N_i}\rlap. \\};
	\arrows (11-) edge node[above]{$\phi$} node[below]{$\simeq$} (-12) (11) edge node[left]{$\psi$} (21) (12-) edge node[above]{$\zeta$} (-13)
	(21-) edge node[below]{$i^\ast(\nu_1)$} node[above]{$\simeq$} (-22) (22-) edge[<-] node[below]{$\simeq$} (-23) (13) edge node[right]{$\nu_1$} node[left]{$\simeq$} (23);
\end{tikzmath}
By the naturality of the isomorphisms~\eqref{eqn:thomses}, it will suffice to verify the commutativity of the following diagram of locally free sheaves:
\begin{tikzmath}
	\diagram{\scr N_i & i^\ast(\Omega_X) & \Omega_{X^f} \\
	\scr N_i & i^\ast(\scr N_\delta) & \scr N_\theta \\
	\scr N_i & \scr N_i & 0\rlap. \\};
	\arrows
	(11-) edge[c->] (-12) (12-) edge[->>] (-13)
	(31-) edge[c->] node[above]{$\id$} (-32) (32-) edge[->>] (-33)
	(11) edge[-,vshift=1pt] (21) edge[-,vshift=-1pt] (21)
	(21) edge node[left]{$\phi$} (31) (12) edge node[right]{$i^\ast(\nu_1)$} node[left]{$\simeq$} (22) (22) edge node[right]{$\psi$} (32)
	(13) edge node[right]{$\nu_1$} node[left]{$\simeq$} (23) (23) edge node[right]{$\zeta$} (33);
\end{tikzmath}
This can be checked on sections as follows. Let $[x]$ be a section of $\scr N_i$, represented by a section $x$ of $\scr O_X$ vanishing on $X^f$. Its image in $i^*(\Omega_X)$ is $i^*(dx)$. By the definitions of $\nu_1$ and $\psi$, we have \[\psi(i^*(\nu_1(dx)))=\psi(i^*[1\tens x-x\tens 1])=[x-x\circ f]=\phi([x]),\] as desired.
This concludes the proof of Theorem~\ref{thm:lefschetz}.

\newpage
\section{The Euler characteristic of separable field extensions}
\label{sec:fields}

In this section we prove Theorem~\ref{thm:trace}. When $L=k$, the statement of Theorem~\ref{thm:trace} reduces to the following lemma:

\begin{lemma}\label{lem:det}
	Let $V$ be a finite-dimensional vector space over $k$ and let $\phi$ be a linear automorphism of $V$. Then $\J\phi=\J{\det(\phi)}$ in $\End(\1_k)$.
\end{lemma}

\begin{proof}
	Recall from \S\ref{sec:review} that $\J\ph$ factors through a group homomorphism $K_1(k)\to\Aut(\1_k)$. The lemma then follows from the fact that the determinant induces an isomorphism $K_1(k)\simeq k^\times$.
\end{proof}

In view of Lemma~\ref{lem:det}, the following proposition completes the proof of Theorem~\ref{thm:trace}:

\begin{proposition}\label{prop:transfers}
	Let $k\subset L$ be a finite separable field extension. For any $\omega\in\End(\1_L)\simeq\GW(L)$,
	\[\int_L\omega\d\chi=\Tr_{L/k}(\omega).\]
	\end{proposition}

\begin{proof}
	Combine Lemmas~\ref{lem:Scharlau}, \ref{lem:geomtransfer}, and \ref{lem:geomScharlau}.
\end{proof}

Recall that, if $p\colon X\to B$ is étale, there are canonical isomorphisms $p^!\simeq p^\ast$ and $p_!\simeq p_\sharp$. If moreover $p$ is finite, we therefore have a canonical isomorphism $p_*\simeq p_\sharp$.

\begin{lemma}\label{lem:etale}
	Let $p\colon X\to B$ be a finite étale morphism and let $\omega\in\End(\1_X)$. Then $\int_X\omega\d\chi\in\End(\1_B)$ is the composition
	\[\1_B\stackrel\eta\to p_*\1_X\simeq p_\sharp\1_X\xrightarrow{p_\sharp\omega}p_\sharp \1_X\stackrel\epsilon\to \1_B.\]
\end{lemma}

\begin{proof}
	By Proposition~\ref{prop:trace}, $\int_X\omega\d\chi$ is the composition
	\begin{tikzmath}
		\def\colsep{1.5em}
		\diagram{\id & p_* p^* & p_* \pi_{1!}\delta_!\delta^!\pi_2^! p^* & p_* \pi_{1!}\pi_2^! p^* & & & \\
		& & & p_!\pi_{2*}\pi_1^*p^! & p_!\pi_{2*}\delta_*\delta^*\pi_1^*p^! & p_!p^! & \id\rlap, \\};
		\arrows (11-) edge node[above]{$\eta$} (-12) (12-) edge node[above]{$\simeq$} (-13) (13-) edge node[above]{$\epsilon$} (-14)
		(14) edge[<-] node[left]{$\Ex_{!*}^{\vphantom{!}}$} node[right]{$\Ex^{*!}_{\vphantom{!}}$} (24)
		(24-) edge node[above]{$\eta$} (-25) (25-) edge node[above]{$\simeq$} (-26) (26-) edge node[above]{$\epsilon$} (-27);
		\draw[<-,font=\scriptsize] (-24) arc (-270:5:10pt);
	\end{tikzmath}
	where the loop is $\omega$ acting after $p^!$. By naturality, we can move this loop to the next-to-last position $p_!p^!$. It then remains to prove that the composition $p_*p^*\to p_!p^!$ (without the loop) is the canonical isomorphism. The morphisms $p$, $\delta$, $\pi_1$, and $\pi_2$ are all finite étale, so we can replace everywhere upper shrieks by upper stars and lower stars by lower shrieks. This operation transforms the exchange isomorphisms $\Ex_{!*}$ and $\Ex^{*!}$ into the connection isomorphisms $\Ex_{!!}$ and $\Ex^{**}$, and we must then prove that the following composition is the identity:
	\[p_!p^*\simeq p_!\pi_{1!}\delta_!\delta^*\pi_{2}^*p^*\to p_!\pi_{1!}\pi_2^*p^*\stackrel c\simeq p_!\pi_{2!}\pi_1^*p^*\to p_!\pi_{2!}\delta_!\delta^*\pi_1^*p^*\simeq p_!p^*.\]
	Using the coherence of the connection isomorphisms, we are reduced to proving that the composition
	\[\delta_!\delta^*\simeq \delta_!\delta^!\stackrel\epsilon\to \id\stackrel\eta\to \delta_*\delta^*\simeq\delta_!\delta^*\]
	is the identity. This is clear since $\delta$ is an open and closed immersion.
	\end{proof}

Fix a base field $k$.
Recall that there is a canonical isomorphism $S^{\A^1}\simeq \P^1_k/\infty$ in $\H_*(k)$ given by the zig-zag
\[\A^1_k/(\A^1_k\minus 0)\to\P^1_k/(\P^1_k\minus 0)\from\P^1_k/\infty.\]

\begin{lemma}\label{lem:identity}
	Let $a\colon\Spec k\into\A^1_k$ be a rational point. Then the composition
	\[S^{\A^1}\simeq \P^1_k/\infty\to\frac{\P^1_k}{\P^1_k\minus a}\simeq S^{\scr N_a}\simeq S^{\A^1}\]
	is the identity in $\H_*(k)$, where the last isomorphism is induced by the trivialization $\scr O_k\simeq\scr N_a$, $1\mapsto t-a$.
\end{lemma}

\begin{proof}
	Suppose first that $a=0$. We must then show that the composition
	\[\frac{\A^1_k}{\A^1_k\minus 0}\to \frac{\P^1_k}{\P^1_k\minus 0}\simeq S^{\scr N_0}\simeq S^{\A^1}\]
	is the identity, which follows from \cite[Lemma 2.2]{VV}. The general case is easily reduced to the case $a=0$ by noting that the map
	\[\P^1/\infty\to \P^1/\infty,\quad [x:y]\mapsto [x+ay:y],\]
	is $\A^1$-homotopic to the identity. 
\end{proof}

\begin{lemma}\label{lem:collapse}
	Let $L$ be a finite separable extension of $k$, $p\colon \Spec L\to\Spec k$ the corresponding morphism of schemes, and $a\colon\Spec L\into \A^1_k$ a closed immersion with minimal polynomial $f\in k[t]$.
	 Then the map $\1_k\stackrel\eta\to p_\ast \1_L\simeq p_\sharp \1_L$ is the $\A^1$-desuspension of the composition
	\[\Sigma^{\A^1}(\Spec k)_+\simeq \P^1_k/\infty\to \frac{\P^1_k}{\P^1_k\minus a}\simeq \Th_{\Spec L}(\scr N_a)
	\simeq \Sigma^{\A^1}(\Spec L)_+,\]
	where $\scr N_a$ is trivialized via $f/f'(a)$.
\end{lemma}

\begin{proof}
	Denote by $\zeta\colon S^{\A^1}\to \P^1_k/(\P^1_k\minus a)$ the first part of the given composition.
	The immersion $a\colon \Spec L\into\P^1_k$ and the given trivialization $\scr N_a\simeq\scr O_L$ form a Euclidean embedding of $\Spec L$ (Definition~\ref{def:euclidean}), and the second part of the given composition is the $\A^1$-suspension of the isomorphism~\eqref{eqn:dual}
	constructed from this Euclidean embedding. By Proposition~\ref{prop:eta}, it therefore suffices to show that the composition
	\begin{equation}\label{eqn:whatwewant}
		\Sigma^{\A^1}(\Spec L)_+\xrightarrow{\zeta\wedge\id}\frac{\P^1_L}{\P^1_L\minus a_L}\stackrel h\to\Sigma^{\A^1}(\Spec L)_+\stackrel p\to S^{\A^1}
	\end{equation}
	is equal to $\Sigma^{\A^1}p_+$ in $\H_*(k)$, where $h$ is the map described in Proposition~\ref{prop:desc} (2). Explicitly, $h$ is the following composition:
	\[\frac{\P^1_L}{\P^1_L\minus a_L}\to\frac{\P^1_L}{\P^1_L\minus \tilde a}\simeq\Th_{\Spec L}(\scr N_{\tilde a})\simeq\Sigma^{\A^1}(\Spec L)_+,\]
	where:
	\begin{itemize}
		\item $a_L\colon\Spec(L\tens_kL)\into \P^1_L$ is the base change of $a$,
		\item $\tilde a=a_L\circ\delta$ is the $L$-point of $\P^1_L$ above $a$,
		\item $\scr N_{\tilde a}$ is trivialized via the isomorphism 
	\begin{equation}\label{eqn:trivNa}
		\scr N_a\simeq a^\ast(\Omega_{\P^1_k})\xrightarrow{\tilde a^\ast (d\hat\pi_1)}\tilde a^\ast(\Omega_{\P^1_L})\simeq \scr N_{\tilde a}
	\end{equation}
	and the given trivialization of $\scr N_a$, where $\hat\pi_1\colon\P^1_L\to\P^1_k$ is the base change of $p$.
	\end{itemize}
	With the identifications
	\[\scr N_a=(f)\tens_{k[t]}L\quad\text{and}\quad \scr N_{\tilde a}=(t-a)\tens_{L[t]}L,\]
	the isomorphism~\eqref{eqn:trivNa} is induced by the inclusion $(f)\subset (t-a)$. If $f(t)=(t-a)g(t)$ in $L[t]$, we have
	\[f(t)\tens \frac 1{f'(a)}=(t-a)g(t)\tens \frac 1{f'(a)}=(t-a)\tens\frac{g(a)}{f'(a)}=(t-a)\tens 1.\]
	Thus, since $\scr N_a$ is trivialized by $f/f'(a)$, $\scr N_{\tilde a}$ is trivialized by the monomial $t-a$.
	The composition
	\[\Sigma^{\A^1}(\Spec L)_+\xrightarrow{\zeta\wedge\id}\frac{\P^1_L}{\P^1_L\minus a_L}\stackrel h\to\Sigma^{\A^1}(\Spec L)_+\]
	is therefore the identity in $\H_*(k)$ by Lemma~\ref{lem:identity} (applied to $\tilde a\colon\Spec L\into \A^1_L$), and hence~\eqref{eqn:whatwewant} is equal to $\Sigma^{\A^1}p_+$, as was to be shown.
\end{proof}

Let $v$ be a finite place of the field of rational functions $k(t)$ with residue field $\kappa(v)$.
As a $k$-vector space, $\kappa(v)$ has a basis $\{1,t,\dotsc,t^{n-1}\}$ where $n=\deg(v)$. We let
\[\tau_v^\mathrm{Sch}\colon \GW(\kappa(v))\to\GW(k)\]
be the Scharlau transfer associated with the $k$-linear map $\kappa(v)\to k$ defined by
\begin{equation*}\label{eqn:scharlau}
	t^i\mapsto\begin{cases}0 & \text{if $0\leq i\leq n-2$,} \\ 1 & \text{if $i=n-1$.}\end{cases}
\end{equation*}
Let also
\[\tau_v^\mathrm{geom}\colon\GW(\kappa(v))\to\GW(k)\]
be the \emph{geometric transfer} defined in \cite[\S4.2]{Morel}.

\begin{lemma}\label{lem:Scharlau}
	Let $v$ be a finite separable place of $k(t)$ with minimal polynomial $f\in k[t]$. Then, for any $\omega\in\GW(\kappa(v))$,
	\[\Tr_{\kappa(v)/k}(\omega)=\tau_v^\mathrm{Sch}(\brac{f'(t)}\omega).\]
\end{lemma}

\begin{proof}
	By \cite[III, \S6, Lemme 2]{Serre}, we have
	\[\Tr_{\kappa(v)/k}\left(\frac{t^i}{f'(t)}\right)=\begin{cases}0 & \text{if $0\leq i\leq n-2$,} \\ 1 & \text{if $i=n-1$.}\end{cases}\]
	This immediately implies the lemma.
	\end{proof}

\begin{lemma}\label{lem:geomtransfer}
	Let $v$ be a finite separable place of $k(t)$ with minimal polynomial $f\in k[t]$. Then, for any $\omega\in\GW(\kappa(v))$,
	\[\int_{\kappa(v)}\omega\d\chi=\tau_v^\mathrm{geom}(\brac{f'(t)}\omega).\]
\end{lemma}

\begin{proof}
	If $a\colon\Spec \kappa(v)\into \P^1_k$ is the closed immersion corresponding to $v$, $\tau_v^\mathrm{geom}$ is the transfer along the same composition as in Lemma~\ref{lem:collapse}, except that the conormal sheaf $\scr N_a$ is trivialized via $f$ (see \cite[\S4.2]{Morel}). The lemma thus follows from Lemmas~\ref{lem:etale} and~\ref{lem:collapse}.
\end{proof}

\begin{lemma}\label{lem:geomScharlau}
	For every finite place $v$ of $k(t)$, $\tau_v^\mathrm{Sch}=\tau_v^\mathrm{geom}$.
\end{lemma}

\begin{proof}
For each place $v$ of $k(t)$, we choose a uniformizer $\pi_v\in\scr O_v$ as follows: if $v$ is finite, let $\pi_v$ be its minimal polynomial, and let $\pi_\infty=-1/t$. By \cite[Theorem 3.15]{Morel}, there is a unique residue homomorphism
\[\del_v\colon K_{\ast +1}^\MW(k(t))\to K_\ast^\MW(\kappa(v))\]
commuting with multiplication by the Hopf element $\eta\in K^\MW_{-1}$ and such that, if $u_1,\dotsc,u_n\in\scr O_v^\times$,
\[\del_v([\pi_v][u_1]\dotso[u_n])=[\bar u_1]\dotso[\bar u_n]\quad\text{and}\quad\del_v([u_1]\dotso[u_n])=0.\]
On the other hand, there are residue homomorphisms
\[\del_v\colon \W(k(t))\to \W(\kappa(v))\]
between Witt groups determined by the formulas
\[\del_v\brac{\pi_v u}=\brac{\bar u}\quad\text{and}\quad\del_v\brac u=0\]
(see \cite[IV, \S1]{HM}). Recalling that $\eta[u]=\brac u-\brac 1$, we see that the following diagram commutes:
\begin{tikzequation}\label{eqn:residues}
	\diagram{K_1^\MW(k(t)) & K_0^\MW(\kappa(v)) \\
	K_1^{\W}(k(t)) & K_0^{\W}(\kappa(v)) \\
	\W(k(t)) & \W(\kappa(v))\rlap. \\};
	\arrows
	(11-) edge node[above]{$\del_v$} (-12)
	(31-) edge node[above]{$\del_v$} (-32)
	(11) edge[->>] (21) (21) edge[c->] node[left]{$\eta$} (31)
	(12) edge[->>] (22) (22) edge[-,vshift=1pt] (32) edge[-,vshift=-1pt] (32);
\end{tikzequation}

By \cite[Theorem 3.24]{Morel}, the map
\[\del=(\del_v)_{v\neq\infty}\colon K_1^\MW(k(t))\to \bigoplus_{v\neq\infty}K_0^\MW(\kappa(v))\]
is surjective, where the sum is taken over all finite places $v$. Given $v$ and $b\in K_0^\MW(\kappa(v))=\GW(\kappa(v))$, choose $\hat b\in K_1^\MW(k(t))$ such that $\del(\hat b)=b$ (in particular, $\del_w(\hat b)=0$ for $w\notin\{v,\infty\}$). By the reciprocity formula for Morel's geometric transfers \cite[(4.8)]{Morel}, we have
\begin{equation}\label{eqn:tau}
	\tau_{v}^\mathrm{geom}(b)=-\del_\infty(\hat b).
\end{equation}
We must therefore show that
\begin{equation}\label{eqn:transfer}
	\tau_v^\mathrm{Sch}(b)=-\del_\infty(\hat b).
\end{equation}
Since $\tau_v^\mathrm{geom}(b)$ and $\tau_v^\mathrm{Sch}(b)$ are both of rank $\deg(v)\cdot\rk(b)$, \eqref{eqn:tau} shows that both sides of~\eqref{eqn:transfer} have the same rank. In view of the cartesian square
\begin{tikzmath}
	\diagram{\GW(k) & \Z \\ \W(k) & \Z/2\rlap, \\};
	\arrows (11-) edge node[above]{$\rk$} (-12) (11) edge (21) (21-) edge (-22) (12) edge (22);
\end{tikzmath}
it remains to prove that~\eqref{eqn:transfer} holds in the Witt group $\W(k)$. By Scharlau's reciprocity theorem for Witt groups \cite[Theorem 4.1]{Scharlau} and~\eqref{eqn:residues}, we have
\[\tau_v^\mathrm{Sch}(b)=\tau_v^\mathrm{Sch}\del_v(\eta\hat b)=-\del_\infty(\eta\hat b)=-\del_\infty(\hat b)\]
in $\W(k)$, as was to be shown. There are two points to be made about the statement of the reciprocity theorem in \cite{Scharlau}. First, the minus sign in front of $\del_\infty$ does not appear in {\it loc.~cit.}, but it appears here because we used the uniformizer $-1/t$ instead of $1/t$ at $\infty$, and we have $\brac{-1/t}=-\brac{1/t}$ in $\W(k(t))$. Second, it is assumed there that $\Char k\neq 2$, but it was observed in \cite[\S2]{GHKS} that, when the Witt group is defined using symmetric bilinear forms instead of quadratic forms, the proof works in arbitrary characteristic.
\end{proof}

\appendix
\newpage
\section{On the purity isomorphism}
\label{app:purity}

In this appendix we achieve two goals:
\begin{itemize}
	\item We show that the purity isomorphism defined by Ayoub \cite[\S1.6]{Ayoub} is the stabilization of the one defined by Morel and Voevodsky \cite[Theorem 2.23]{MV}.
	\item We prove that the purity isomorphism $\Pi\colon s^!p^*\simeq\Sigma^{-\scr N_s}q^*$ is natural in the closed immersion $s$.
	\end{itemize}
The naturality of $\Pi$ plays a central role in the proof of Theorem~\ref{thm:lefschetz} in \S\ref{sec:traces}.

We start by recalling the definition of the Morel–Voevodsky purity zig-zag. Let $C$ be the open subscheme of the blowup of $X\times\A^1$ along $Z\times\{0\}$ whose closed complement is the blowup of $X\times\{0\}$ along $Z\times\{0\}$. We then have canonical isomorphisms
	\[C\times_{\A^1}\{0\}\simeq \bb V(\scr N_s)\quad\text{and}\quad C\times_{\A^1}\{1\}\simeq X\]
	(see \cite[Chapter 5]{Fulton}) and diagrams
	\[
	\begin{tikzpicture}
		\diagram{Z & X \\ & B \\};
		\arrows (11-) edge[c->] node[above]{$s$} (-12) (11) edge (22) (12) edge node[right]{$p$} (22);
	\end{tikzpicture}
	\quad
	\begin{tikzpicture}
		\draw[c->,font=\scriptsize] (0,0) -- node[above]{$i$} (1,0);
	\end{tikzpicture}
	\quad
	\begin{tikzpicture}
		\diagram{Z\times\A^1 & C \\ & B\times\A^1 \\};
		\arrows (11-) edge[c->] node[above]{$\hat s$} (-12) (11) edge (22) (12) edge node[right]{$\hat p$} (22);
	\end{tikzpicture}
	\quad
	\begin{tikzpicture}
		\draw[left hook->,font=\scriptsize] (1,0) -- node[above]{$i_0$} (0,0);
	\end{tikzpicture}
	\quad
	\begin{tikzpicture}
		\diagram{Z & \bb V(\scr N_s) \\ & B \\};
		\arrows (11-) edge[c->] node[above]{$s_0$} (-12) (11) edge (22) (12) edge node[right]{$p_0$} (22);
	\end{tikzpicture}
	\]
	where $i$ (\resp{} $i_0$) is the inclusion of the fiber over $1$ (\resp{} over $0$). Note that $\Sigma^{-\scr N_s}s^\ast p^\ast\simeq s_0^!p_0^\ast$ by definition of $\Sigma^{-\scr N_s}$.
	Denote by $r\colon B\times\A^1\to B$ and $r\colon Z\times\A^1\to Z$ the projections. Since $i$ is a section of $r$, there is a transformation $r_\ast\to i^\ast$ given by
	\begin{equation}\label{eqn:ri}
		r_\ast\stackrel\eta\to r_\ast i_\ast i^\ast \simeq i^\ast.
	\end{equation}
	Let $\Pi_1$ be the composition
	\begin{equation}\label{eqn:Pi}
	r_\ast \hat s^! \hat p^\ast r^\ast \xrightarrow{\eqref{eqn:ri}} i^\ast\hat s^! \hat p^\ast r^\ast\xrightarrow{\Ex^{\ast!}} s^!i^\ast \hat p^\ast r^\ast\simeq s^!p^\ast,
	\end{equation}
	and let $\Pi_0\colon r_\ast \hat s^! \hat p^\ast r^\ast\to s_0^!p_0^\ast$ be the analogous composition with $i$ replaced by $i_0$, so that we have a zig-zag
	\[s^!p^*\xleftarrow{\Pi_1} r_\ast \hat s^! \hat p^\ast r^\ast\xrightarrow{\Pi_0} s_0^!p_0^*.\]
	
\begin{proposition}\label{prop:MV=Ayoub}
	The transformations $\Pi_1$ and $\Pi_0$ are isomorphisms and the composition $\Pi_0\Pi_1^{-1}$ coincides with the purity isomorphism $\Pi$.
\end{proposition}

\begin{proof}
	We will show that both maps in~\eqref{eqn:Pi} are isomorphisms. Consider the diagram
	\begin{tikzmath}
		\diagram{r_\ast \hat s^! \hat p^\ast r^\ast & i^\ast\hat s^! \hat p^\ast r^\ast & s^!p^\ast \\
		r_\ast \Sigma^{-\scr N_{\hat s}} \hat s^\ast\hat p^\ast r^\ast & i^\ast \Sigma^{-\scr N_{\hat s}} \hat s^\ast\hat p^\ast r^\ast & \Sigma^{-\scr N_s}s^\ast p^\ast\rlap. \\};
		\arrows (11-) edge node[above]{\eqref{eqn:ri}} (-12) (12-) edge node[above]{$\Ex^{\ast !}$} (-13)
		(21-) edge node[above]{\eqref{eqn:ri}} (-22) (22-) edge node[above]{$\simeq$} (-23)
		(11) edge node[left]{$\Pi$} node[right]{$\simeq$} (21) (12) edge node[left]{$\Pi$} node[right]{$\simeq$} (22) (13) edge node[left]{$\Pi$} node[right]{$\simeq$} (23);
	\end{tikzmath}
	The first square commutes by naturality of the transformation~\eqref{eqn:ri} and the second square commutes by \cite[Corollaire 1.6.23]{Ayoub}. Moreover, the transformation at the bottom left is an isomorphism because $\scr N_{\hat s}\simeq r^\ast(\scr N_s)$ and $\eqref{eqn:ri}r^\ast$ is an isomorphism. Using that $\eta\colon\id\to r_\ast r^\ast$ is an isomorphism, we see that the lower row does not change if we replace $i$ by $i_0$. Together with the analogous diagram for $\Pi_0$, we therefore obtain a commutative square
	\begin{tikzmath}
		\diagram{s^!p^\ast & s_0^!p_0^\ast \\ s_0^!p_0^\ast & s_0^!p_0^\ast\rlap.\\};
		\arrows (11-) edge node[above]{$\Pi_0\Pi_1^{-1}$} (-12) (11) edge node[left]{$\Pi$} (21)
		(12) edge node[right]{$\Pi$} (22) (21-) edge node[below]{$\id$} (-22);
	\end{tikzmath}
	But the right-hand purity isomorphism $\Pi\colon s_0^!p_0^\ast\simeq s_0^!p_0^\ast$ is the identity by \cite[Proposition 1.6.28]{Ayoub}, and hence $\Pi=\Pi_0\Pi_1^{-1}$, as claimed.
\end{proof}

\begin{proposition}[Purity is natural]
	\label{prop:purity}
	Suppose given a cartesian square
	\begin{tikzmath}
		\diagram{W & Y \\ Z & X \\};
		\arrows (11-) edge[c->] node[above]{$t$} (-12)
		(21-) edge[c->] node[above]{$s$} (-22)
		(11) edge node[left]{$g$} (21) (12) edge node[right]{$f$} (22);
	\end{tikzmath}
	in $\Sm_B$ where $s$ and $t$ are closed immersions, and let $p\colon X\to B$ be the structure map. Then the induced map $\psi\colon g^\ast(\scr N_s)\to\scr N_t$ is an epimorphism and the diagrams
	\[
	\begin{tikzpicture}
		\diagram{g^\ast s^! p^\ast & g^\ast\Sigma^{-\scr N_s} s^\ast p^\ast & \Sigma^{-g^\ast(\scr N_s)}g^\ast s^\ast p^\ast \\
		t^! f^\ast p^\ast & \Sigma^{-\scr N_t}t^\ast f^\ast p^\ast & \Sigma^{-\scr N_t} g^\ast s^\ast p^\ast \\};
		\arrows
		(11) edge node[left]{$\Ex^{\ast !}$} (21)
		(13) edge node[right]{$\Sigma^{-\psi}$} (23)
		(11-) edge node[above]{$\Pi$} node[below]{$\simeq$} (-12)
		(21-) edge node[above]{$\Pi$} node[below]{$\simeq$} (-22)
		(12-) edge node[below]{$\simeq$} (-13)
		(22-) edge node[below]{$\simeq$} (-23)
		;
	\end{tikzpicture}
	\quad
	\begin{tikzpicture}
		\diagram{g^! s^* p^! & g^!\Sigma^{\scr N_s} s^! p^! & \Sigma^{g^*(\scr N_s)}g^! s^! p^! \\
		t^* f^! p^! & \Sigma^{\scr N_t}t^! f^! p^! & \Sigma^{\scr N_t} g^! s^! p^! \\};
		\arrows
		(11) edge[<-] node[left]{$\Ex^{\ast !}$} (21)
		(13) edge[<-] node[right]{$\Sigma^{\psi}$} (23)
		(11-) edge node[above]{$\Pi$} node[below]{$\simeq$} (-12)
		(21-) edge node[above]{$\Pi$} node[below]{$\simeq$} (-22)
		(12-) edge node[below]{$\simeq$} (-13)
		(22-) edge node[below]{$\simeq$} (-23)
		;
	\end{tikzpicture}
	\]
	are commutative.
\end{proposition}

\begin{remark}
	When $f$ is smooth (in which case $\psi$ is an isomorphism), Proposition~\ref{prop:purity} is exactly \cite[Proposition 1.6.20]{Ayoub}, but in~\S\ref{sec:traces} we need the proposition for $f$ a closed immersion.
\end{remark}

\begin{proof}
	Let $\scr I\subset\scr O_X$ be the defining ideal of $s$ and $\scr J\subset\scr O_Y$ that of $t$. The morphism $\psi$ is then the composition
	\[g^\ast(\scr N_s)\simeq g^\ast s^\ast(\scr I)\simeq t^\ast f^\ast(\scr I)\to t^\ast(\scr J)\simeq\scr N_t.\]
	Because the square is cartesian, $\scr J$ is exactly the image of $f^\ast(\scr I)\to\scr O_Y$, and since $t^\ast$ is right exact, $\psi$ is an epimorphism.
	
	We will only prove the commutativity of the first diagram; the commutativity of the second diagram is checked by a dual argument.
	Let $D$ be the open subscheme of the blowup of $Y\times\A^1$ along $W\times\{0\}$ whose complement is the proper transform of $Y\times\{0\}$. The given cartesian square then induces cartesian squares
	\[
	\begin{tikzpicture}
		\diagram{W & Y \\ Z & X \\};
		\arrows (11-) edge[c->] node[above]{$t$} (-12)
		(21-) edge[c->] node[above]{$s$} (-22)
		(11) edge node[left]{$g$} (21) (12) edge node[right]{$f$} (22);
	\end{tikzpicture}
	\quad
	\begin{tikzpicture}
		\draw[c->,font=\scriptsize] (0,0) -- node[above]{$i$} (1,0);
	\end{tikzpicture}
	\quad
	\begin{tikzpicture}
		\diagram{W\times\A^1 & D \\ Z\times\A^1 & C \\};
		\arrows (11-) edge[c->] node[above]{$\hat t$} (-12)
		(21-) edge[c->] node[above]{$\hat s$} (-22)
		(11) edge node[left]{$\hat g$} (21) (12) edge node[right]{$\hat f$} (22);
	\end{tikzpicture}
	\quad
	\begin{tikzpicture}
		\draw[left hook->,font=\scriptsize] (1,0) -- node[above]{$i_0$} (0,0);
	\end{tikzpicture}
	\quad
	\begin{tikzpicture}
		\diagram{W & \bb V(\scr N_t) \\ Z & \bb V(\scr N_s) \\};
		\arrows (11-) edge[c->] node[above]{$t_0$} (-12)
		(21-) edge[c->] node[above]{$s_0$} (-22)
		(11) edge node[left]{$g_0$} (21) (12) edge node[right]{$f_0$} (22);
	\end{tikzpicture}
	\]
	where $f_0=\bb V(\psi)$ and $g_0=g$. By Lemma~\ref{lem:thommap}, the transformation $\Sigma^{-\psi}\colon \Sigma^{-g^\ast\scr N_s} g^\ast s^\ast p^\ast\to \Sigma^{-\scr N_t} g^\ast s^\ast p^\ast$ can be identified with the exchange transformation $\Ex^{\ast !}\colon g_0^\ast s_0^!p_0^\ast\to t_0^! f_0^\ast p_0^\ast$. 
	By replacing both occurrences of $\Pi$ by $\Pi_0\Pi_1^{-1}$ (Proposition~\ref{prop:MV=Ayoub}) and completing the resulting diagram with exchange transfomations of the form $\Ex^{\ast !}$, we are reduced to proving the commutativity of the rectangle
	\begin{tikzmath}
		\diagram{i^\ast \hat g^\ast\hat s^! & g^\ast i^\ast\hat s^! & g^\ast s^! i^\ast \\ i^\ast\hat t^!\hat f^\ast & t^!i^\ast\hat f^\ast & t^!f^\ast i^\ast \\};
		\arrows (11-) edge node[above]{$\simeq$} (-12) (21-) edge node[above]{$\Ex^{\ast !}$} (-22)
		(12-) edge node[above]{$\Ex^{\ast !}$} (-13) (11) edge node[left]{$\Ex^{\ast !}$} (21) (22-) edge node[above]{$\simeq$} (-23) (13) edge node[right]{$\Ex^{\ast !}$} (23);
	\end{tikzmath}
	and of the analogous rectangle with $i$ replaced by $i_0$.
	By formal properties of exchange transformations \cite[Définition 1.2.1]{Ayoub}, both compositions in this rectangle are equal to
	\[i^*\hat g^* s^!\stackrel c\simeq (\hat gi)^*s^!= (ig)^* s^!\xrightarrow{\Ex^{*!}} t^!(if)^*\stackrel c\simeq t^!f^*i^*.\qedhere\]
	\end{proof}

\newpage
\section{Coherence lemmas}
\label{app:lemmas}

\begin{lemma}\label{lem:bc}
	Let
	\begin{tikzmath}
		\diagram{X\times_BX & X \\ X & B \\};
		\arrows (11-) edge node[above]{$\pi_2$} (-12) (21-) edge node[below]{$p$} (-22) (11) edge node[left]{$\pi_1$} (21) (12) edge node[right]{$p$} (22);
	\end{tikzmath}
	be a cartesian square of schemes, and let $\delta\colon X\into X\times_BX$ be the diagonal. Then:
	\begin{enumerate}
		\item The counit $\epsilon\colon p^*p_*\to\id$ coincides with the composition
		\[p^*p_*\xrightarrow{\Ex^*_*} \pi_{2*}\pi_1^*\stackrel\eta\to\pi_{2*}\delta_*\delta^*\pi_1^* \simeq\id.\]
		\item If $p$ is separated of finite type, then the unit $\eta\colon\id\to p^!p_!$ coincides with the composition
		\[\id\simeq \pi_{1!}\delta_!\delta^!\pi_2^!\stackrel \epsilon\to \pi_{1!}\pi_2^!\xrightarrow{\Ex_!^!} p_!p^!.\]
	\end{enumerate}
\end{lemma}

\begin{proof}
	The diagram
	\begin{tikzmath}
		\diagram{
		p^*p_* & p^*p_*\pi_{1*}\pi_1^* & p^*p_*\pi_{2*}\pi_1^* & \pi_{2*}\pi_1^* \\
		 & p^*p_*\pi_{1*}\delta_*\delta^*\pi_1^* & p^*p_*\pi_{2*}\delta_*\delta^*\pi_1^* & \pi_{2*}\delta_*\delta^*\pi_1^* \\
		 & & p^*p_*\id_*\id^* & \id_*\id^* \\};
		\arrows (11-) edge node[above]{$\eta$} (-12) (12-) edge node[above]{$c$} (-13) (13-) edge node[above]{$\epsilon$} (-14)
		(22-) edge node[above]{$c$} (-23) (23-) edge node[above]{$\epsilon$} (-24)
		(22) edge node[below left]{$c$} (33) (33-) edge node[above]{$\epsilon$} (-34)
		(12) edge node[left]{$\eta$} (22) (13) edge node[left]{$\eta$} (23) (14) edge node[left]{$\eta$} (24)
		(23) edge node[left]{$c$} (33) (24) edge node[left]{$c$} (34);
		\draw[->,font=\scriptsize] (11) to[out=290,in=180] node[below left]{$\eta$} (33);
	\end{tikzmath}
	is clearly commutative. Comparing the two outer compositions proves (1). The proof of (2) is identical.
\end{proof}

\begin{lemma}\label{lem:symmetry}
	Let
	\begin{tikzmath}
		\diagram{\bullet & \bullet \\ \bullet & \bullet \\};
		\arrows (11-) edge node[above]{$g$} (-12) (11) edge node[left]{$q$} (21) (21-) edge node[below]{$f$} (-22) (12) edge node[right]{$p$} (22);
	\end{tikzmath}
	be a cartesian square of schemes where $p$ is proper and $f$ is smooth. Then the following diagram commutes:
	\begin{tikzmath}
		\def\colsep{3.5em}
		\diagram{p_!p^*\1_B\wedge f_*f^!\1_B & p_!p^*f_*f^!\1_B & p_!g_{*}q^*f^!\1_B \\
		f_*f^!\1_B\wedge p_!p^*\1_B & f_*f^!p_!p^*\1_B & f_*q_{!}g^!p^*\1_B\rlap. \\};
		\arrows (11-) edge[<-] node[above]{$\Pr^*_!$} (-12) (21-) edge node[above]{$\Pr^{*!}\Pr^*_*$} (-22)
		(12-) edge node[above]{$\Ex^*_*$} (-13) (-23) edge node[above]{$\Ex_!^!$} (22-)
		(11) edge node[left]{$\tau$} (21) (13) edge node[left]{$\Ex_{!*}^{\vphantom{!}}$} node[right]{$\Ex^{*!}_{\vphantom{!}}$} (23);
	\end{tikzmath}
\end{lemma}

\begin{proof}
	Since $p$ is proper and $f$ is smooth, we can eliminate the shrieks to obtain the equivalent rectangle:
	\begin{tikzmath}
		\diagram{p_*p^*\1_B\wedge f_*\Sigma^{\Omega_f} f^*\1_B & p_*p^*f_*\Sigma^{\Omega_f} f^*\1_B & p_*g_*\Sigma^{\Omega_g} q^*f^*\1_B \\
		f_*\Sigma^{\Omega_f} f^*\1_B\wedge p_*p^*\1_B & f_*\Sigma^{\Omega_f} f^*p_*p^*\1_B & f_* q_*\Sigma^{\Omega_g}g^*p^*\1_B\rlap. \\};
		\arrows (11-) edge node[above]{$\Pr^*_*$} (-12) (21-) edge node[above]{$\Pr^*_*$} (-22)
		(12-) edge node[above]{$\Ex^*_*$} (-13) (-23) edge[<-] node[above]{$\Ex_* ^*$} (22-)
		(11) edge node[left]{$\tau$} (21) (13) edge node[left]{$\Ex_{* *}^{\vphantom{!}}$} node[right]{$\Ex^{**}_{\vphantom{!}}$} (23);
	\end{tikzmath}
	This is now a special case of the following diagram, for $E=p^*\1_B$ and $F=\Sigma^{\Omega_f}f^*\1_B$:
	\begin{tikzmath}
		\diagram{p_*E\wedge f_* F & p_*(E\wedge p^* f_* F) & p_*(E\wedge g_* q^*F) & p_*g_*(g^*E\wedge q^*F) \\
		f_* F\wedge p_*E & f_* (F\wedge f^*p_*E) & f_*(F\wedge q_*g^*E) & f_* q_*(q^*F\wedge g^*E)\rlap. \\};
		\arrows (11-) edge node[above]{$\Pr^*_*$} (-12) (21-) edge node[above]{$\Pr^*_*$} (-22)
		(12-) edge node[above]{$\Ex^*_*$} (-13) (-23) edge[<-] node[above]{$\Ex_* ^*$} (22-)
		(13-) edge node[above]{$\Pr^*_*$} (-14) (23-) edge node[above]{$\Pr^*_*$} (-24)
		(11) edge node[left]{$\tau$} (21) (14) edge node[left]{$c$} node[right]{$\tau$} (24);
	\end{tikzmath}
	We will prove that this more general rectangle is commutative. Using the adjunction $(fq)^*\dashv (fq)_*$ and the fact that the functors $(\ph)^*$ are symmetric monoidal, it is equivalent to prove the commutativity of the following rectangle:
	\begin{tikzmath}
		\diagram{g^*p^*p_*E\wedge g^*p^*f_* F & g^*E\wedge g^*p^* f_* F & g^*E\wedge g^*g_* q^*F & g^*E\wedge q^*F \\
		q^*f^*f_* F\wedge q^*f^*p_*E & q^* F\wedge q^*f^*p_*E & q^* F\wedge q^*q_*g^*E & q^* F\wedge g^*E\rlap. \\};
		\arrows (11-) edge node[above]{$\epsilon\wedge\id$} (-12) (21-) edge node[above]{$\epsilon\wedge\id$} (-22)
		(12-) edge node[above]{$\Ex^*_*$} (-13) (22-) edge node[above]{$\Ex_* ^*$} (-23)
		(13-) edge node[above]{$\id\wedge\epsilon$} (-14) (23-) edge node[above]{$\id\wedge\epsilon$} (-24)
		(11) edge node[left]{$\tau\circ(c\wedge c)$} (21) (14) edge node[right]{$\tau$} (24);
	\end{tikzmath}
	The compositions in this rectangle are now of the form $\tau\circ(\phi\wedge\psi)$ and $(\psi'\wedge\phi')\circ\tau$, and hence we need only check that $\phi=\phi'$ and $\psi= \psi'$, \ie, that the squares
	\[
	\begin{tikzpicture}
		\diagram{g^*p^*f_* & g^*g_*q^* \\ q^*f^*f_* & q^* \\};
		\arrows (11-) edge node[above]{$\Ex_*^*$} (-12) (11) edge node[left]{$c$} (21) (21-) edge node[below]{$\epsilon$} (-22) (12) edge node[right]{$\epsilon$} (22);
	\end{tikzpicture}
	\qquad
	\begin{tikzpicture}
		\diagram{q^*f^*p_* & q^*q_*g^* \\ g^*p^*p_* & g^* \\};
		\arrows (11-) edge node[above]{$\Ex_*^*$} (-12) (11) edge node[left]{$c$} (21) (21-) edge node[below]{$\epsilon$} (-22) (12) edge node[right]{$\epsilon$} (22);
	\end{tikzpicture}
	\]
	are commutative. This follows from \cite[Proposition 1.2.5]{Ayoub}.
	\end{proof}

\begin{lemma}\label{lem:prism}
	Let
	\begin{tikzmath}
		\diagram{\bullet & \bullet & \bullet \\ \bullet & \bullet & \bullet \\};
		\arrows (11-) edge node[above]{$g$} (-12) (11) edge node[left]{$r''$} (21) (21-) edge node[below]{$f$} (-22) (12-) edge node[above]{$q$} (-13) (12) edge node[left]{$r'$} (22) (22-) edge node[below]{$p$} (-23) (13) edge node[right]{$r$} (23);
	\end{tikzmath}
	be cartesian squares in which all maps are separated of finite type. Then the following rectangle commutes:
	\begin{tikzmath}
		\diagram{(pf)_!r''_!(qg)^! & p_!f_!r''_!g^{!}q^! & p_!r'_!g_! g^{!}q^{!} & p_!r'_!q^{!} \\
		(pf)_!(pf)^!r_! & p_!f_!f^!p^!r_! & & p_!p^!r_!\rlap. \\};
		\arrows (11-) edge node[above]{$c$} (-12) (12-) edge node[above]{$\Ex_{!!}$} (-13) (13-) edge node[above]{$\epsilon$} (-14)
		(11) edge node[left]{$\Ex_!^!$} (21) (21-) edge node[below]{$c$} (-22) (22-) edge node[below]{$\epsilon$} (-24) (14) edge node[right]{$\Ex_!^!$} (24);
	\end{tikzmath}
\end{lemma}

\begin{proof}
	We break up this rectangle as follows:
	\begin{tikzmath}
		\diagram{(pf)_!r''_!(qg)^! & p_!f_!r''_!g^{!}q^! & p_!r'_!g_! g^{!}q^{!} \\
		& p_!f_!f^!r'_!q^{!} & p_!r'_!q^{!} \\
		(pf)_!(pf)^!r_! & p_!f_!f^!p^!r_! & p_!p^!r_!\rlap. \\};
		\arrows (11-) edge node[above]{$c$} (-12) (12-) edge node[above]{$\Ex_{!!}$} (-13) (13) edge node[right]{$\epsilon$} (23)
		(12) edge node[left]{$\Ex_!^!$} (22) (22) edge node[left]{$\Ex_!^!$} (32)
		(22-) edge node[above]{$\epsilon$} (-23)
		(11) edge node[left]{$\Ex_!^!$} (31) (31-) edge node[below]{$c$} (-32) (32-) edge node[below]{$\epsilon$} (-33)
		(23) edge node[right]{$\Ex_!^!$} (33);
	\end{tikzmath}
	The left rectangle commutes by the compatibility of exchange transformations with the composition of cartesian squares \cite[Définition 1.2.1]{Ayoub}, the top square commutes by \cite[Proposition 1.2.5]{Ayoub}, and the bottom square commutes by naturality of $\epsilon$.
\end{proof}

\begin{lemma}\label{lem:bottomrow}
	Let $\gamma,\delta\colon X\into Y$ be a pair of closed immersions with a common retraction $\pi\colon Y\to X$, and let
	\begin{tikzmath}
		\diagram{Z & X \\ X & Y \\};
		\arrows (11-) edge node[above]{$i$} (-12) (11) edge node[left]{$i$} (21) (12) edge node[right]{$\delta$} (22) (21-) edge node[below]{$\gamma$} (-22);
	\end{tikzmath}
	be a cartesian square. Then the composition
	\[i_*i^!\xrightarrow{\Ex_*^!} \gamma^!\delta_* \simeq \pi_*\gamma_*\gamma^!\delta_*\stackrel\epsilon\to \pi_*\delta_*\simeq\id\]
	is equal to the counit $\epsilon\colon i_*i^!\to\id$.
	\end{lemma}

\begin{proof}
	Consider the commutative diagram
	\begin{tikzmath}
		\diagram{i_*i^! & \gamma^!\gamma_*i_*i^! & \gamma^!\delta_*i_*i^! & \gamma^!\delta_* \\
		\pi_{*}\gamma_*i_*i^! & \pi_{*}\gamma_*\gamma^!\gamma_*i_*i^! & \pi_{*}\gamma_*\gamma^!\delta_*i_*i^! & \pi_{*}\gamma_* \gamma^!\delta_* \\
		 & \pi_{*}\gamma_*i_*i^! & \pi_{*}\delta_*i_*i^! & \pi_{*}\delta_*\rlap{${}\simeq \id$,} \\};
		 \arrows
		 (11-) edge node[above]{$\eta$} (-12) (12-) edge node[above]{$c$} (-13) (13-) edge node[above]{$\epsilon$} (-14)
		 (21-) edge node[above]{$\eta$} (-22) (22-) edge node[above]{$c$} (-23) (23-) edge node[above]{$\epsilon$} (-24)
		 (32-) edge node[above]{$c$} (-33) (33-) edge node[above]{$\epsilon$} (-34)
		 (11) edge node[left]{$c$} (21) (12) edge node[left]{$c$} (22) (13) edge node[left]{$c$} (23) (14) edge node[left]{$c$} (24)
		 (21) edge node[below left]{$\id$} (32) (22) edge node[left]{$\epsilon$} (32) (23) edge node[left]{$\epsilon$} (33) (24) edge node[left]{$\epsilon$} (34);
		 ;
	\end{tikzmath}
	in which the upper composition is the given one.
	By coherence of the connection isomorphisms, the lower composition is the counit $\epsilon\colon i_*i^!\to\id$, which proves the lemma.
\end{proof}

\begin{lemma}\label{lem:sigma}
	Let $i$ be a closed immersion. Then the triangle
	\begin{tikzmath}
		\diagram{i^! & i^* \\ & i^!i_*i^* \\};
		\arrows (11-) edge node[above]{$\sigma$} (-12) (11) edge node[below left]{$\eta$} (22) (12) edge node[left]{$\simeq$} node[right]{$\eta$} (22);
	\end{tikzmath}
	is commutative.
\end{lemma}

\begin{proof}
	This is simply a matter of unwinding the definitions.
	Recall that $\sigma$ is $\Ex^{*!}\colon \id^*i^!\to \id^!i^*$. By definition of $\Ex^{*!}$, this is the composition
	\[\id^*i^!\stackrel\eta\to \id^*i^!i_*i^*\xleftarrow{\Ex^!_*} \id^*\id_*\id^! i^*\stackrel\epsilon\to \id^!i^*.\]
	By construction, $\Ex^!_*\colon \id_*\id^!\to i^!i_*$ is the mate of $\Ex^*_*\colon i^*i_*\to\id_*\id^*$. Finally, by definition of $\Ex^*_*$, the latter is the counit $\epsilon\colon i^*i_*\to\id$, whose mate is $\eta\colon\id\to i^!i_*$.
\end{proof}

\begin{lemma}\label{lem:thommap}
	Let
	\begin{tikzmath}
		\diagram{W & V \\ Y & X \\};
		\arrows (11-) edge node[above]{$g$} (-12) (11) edge node[left]{$q$} (21)
		(12) edge node[right]{$p$} (22) (21-) edge node[below]{$f$} (-22);
	\end{tikzmath}
	be a commutative diagram, where $p$ and $q$ are vector bundles with zero sections $s$ and $t$, and where $g$ induces a monomorphism of vector bundles $\phi\colon W\into f^*V$. Then the following diagrams commute (the second assuming that $f$ is separated of finite type):
	\[
	\begin{tikzpicture}
		\diagram{f^*\Sigma^{-V} & \Sigma^{-f^*V}f^* & \Sigma^{-W}f^* \\
		f^*s^!p^* & t^!g^* p^* & t^! q^* f^*\rlap, \\};
		\arrows (11-) edge node[above]{$\simeq$} (-12) (12-) edge node[above]{$\Sigma^{-\phi}$} (-13)
		(11) edge node[left]{$\simeq$} (21)
		(13) edge node[right]{$\simeq$} (23)
		(21-) edge node[below]{$\Ex^{*!}$} (-22) (22-) edge node[below]{$\simeq$} (-23);
	\end{tikzpicture}
	\quad
	\begin{tikzpicture}
		\diagram{f^!\Sigma^{V} & \Sigma^{f^*V}f^! & \Sigma^{W}f^! \\
		f^!s^*p^! & t^*g^! p^! & t^* q^! f^!\rlap. \\};
		\arrows (11-) edge[<-] node[above]{$\simeq$} (-12) (12-) edge[<-] node[above]{$\Sigma^{\phi}$} (-13)
		(11) edge[<-] node[left]{$\simeq$} (21)
		(13) edge[<-] node[right]{$\simeq$} (23)
		(21-) edge[<-] node[below]{$\Ex^{*!}$} (-22) (22-) edge[<-] node[below]{$\simeq$} (-23);
	\end{tikzpicture}
	\]
\end{lemma}

\begin{proof}
	Let $r\colon f^*V\to Y$ be the pullback of $p$ and let $u$ be the zero section of $r$. Recall that $\Sigma^{-\phi}$ is the composition
	\[u^!r^*\stackrel c\simeq t^!\phi^!r^*\stackrel \sigma\to t^!\phi^*r^*\stackrel c\simeq t^! q^*,\]
	and that $\sigma\colon \phi^!\to \phi^*$ is the exchange transformation $\Ex^{*!}\colon \id^*\phi^!\to \id^!\phi^*$. The commutativity of the first rectangle then follows from the compatibility of the exchange transformation $\Ex^{*!}$ with the composition of the following three cartesian squares:
	\begin{tikzmath}
		\diagram{Y & W & W \\ Y & W & f^*V \\ X & & V\rlap. \\};
		\arrows (11-) edge node[above]{$t$} (-12) (12-) edge[-,vshift=1.2pt] (-13) edge[-,vshift=-1.2pt] (-13)
		(11) edge[-,vshift=1.2pt] (21) edge[-,vshift=-1.2pt] (21)
		(12) edge[-,vshift=1.2pt] (22) edge[-,vshift=-1.2pt] (22)
		(13) edge node[right]{$\phi$} (23) (21-) edge node[below]{$t$} (-22) (22-) edge node[below]{$\phi$} (-23)
		(21) edge node[left]{$f$} (31) (23) edge (33)
		(31-) edge node[above]{$s$} (-33);
	\end{tikzmath}
	The commutativity of the second square is checked in the same way.
\end{proof}

\newpage
\section{Elimination of noetherian hypotheses}
\label{app:qcqs}

In the foundational paper \cite{MV}, unstable motivic homotopy theory is only defined for noetherian schemes of finite Krull dimension. In this appendix we indicate how to properly extend the theory to arbitrary schemes. For simplicity, we will give our definitions using the language of $\infty$-categories \cite{HTT}. We say that a scheme or a morphism of schemes is \emph{coherent} if it is quasi-compact and quasi-separated.

There are two issues that arise when dropping the assumption that schemes are noetherian and finite-dimensional.
The first concerns the definition of the Nisnevich topology. This topology was originally defined in \cite{Nisnevich} using the following pretopology: a family $\{U_i\to X\}_{i\in I}$ is a cover if each $U_i\to X$ is étale and every morphism $\Spec k\to X$ with $k$ a field lifts to $U_i$ for some $i\in I$. 
For noetherian schemes, it was shown in \cite[Proposition 3.1.4]{MV} that this topology is generated by a cd-structure, in the sense of \cite[\S2]{VV:cd}. For coherent schemes that are not noetherian, the pretopology and the cd-structure define different topologies, both finer than the Zariski topology and coarser than the étale topology.
We will define the Nisnevich topology in general by combining the cd-structure and the Zariski topology.
This choice ensures that the ``small'' Nisnevich $\infty$-topos $X_\Nis$ of a scheme $X$
(\ie, the $\infty$-category of Nisnevich sheaves of spaces on étale $X$-schemes)
has good formal properties. For instance:
\begin{enumerate}
	\item if $X$ is coherent, then $X_\Nis$ is coherent and compactly generated by finitely presented étale $X$-schemes;
	\item if $X$ is the limit of a cofiltered diagram of coherent schemes $X_\alpha$ with affine transition maps, then $X_\Nis$ is the limit of the $\infty$-topoi $(X_\alpha)_\Nis$.
\end{enumerate}
Another point in favor of our definition is that algebraic $K$-theory, considered as a presheaf of spaces on coherent schemes, is only known to be a sheaf for our version of the Nisnevich topology.
Note that property (2) determines the $\infty$-topos $X_\Nis$ for $X$ coherent once it has been defined for $X$ noetherian, since any coherent scheme is a cofiltered limit of schemes of finite type over $\Z$ \cite[Appendix C]{Thomason}. The second issue is that the Nisnevich $\infty$-topos of a coherent scheme which is not noetherian and finite-dimensional need not be hypercomplete, \ie, Nisnevich descent for a presheaf of spaces does not imply Nisnevich hyperdescent. We do not want to restrict ourselves to hypercomplete sheaves, since by doing so we might lose properties (1) and (2) as well as the representability of algebraic $K$-theory.

In this appendix, a presheaf is by default a presheaf of \emph{spaces}. If $\scr C$ is a (possibly large) $\infty$-category, we denote by $\PSh(\scr C)$ the $\infty$-category of presheaves on $\scr C$. It will be convenient to work with a weakening of the notion of topology: a \emph{quasi-topology} $\tau$ on an $\infty$-category $\scr C$ assigns to every $X\in\scr C$ a collection $\tau(X)$ of sieves on $X$, called \emph{$\tau$-sieves}, such that, for every $f\colon Y\to X$, $f^*\tau(X)\subset\tau (Y)$. A presheaf $F$ on $\scr C$ is a \emph{$\tau$-sheaf} if, for every $X\in\scr C$ and every $R\in \tau(X)$, the restriction map $\Map(X,F)\to\Map(R,F)$ is an equivalence. We denote by $\Shv_\tau(\scr C)\subset\PSh(\scr C)$ the full subcategory of $\tau$-sheaves. A family of morphisms $\{U_i\to X\}$ in $\scr C$ is called a \emph{$\tau$-cover} if it generates a $\tau$-sieve.

If $\tau$ is a quasi-topology on $\scr C$, we denote by $\bar\tau$ the coarsest topology containing $\tau$. Our first goal is to show that $\Shv_\tau(\scr C)=\Shv_{\bar\tau}(\scr C)$. The following proposition is a generalization of \cite[II, Proposition 2.2]{SGA4-1} to sheaves of spaces; the proof is exactly the same.
 
\begin{proposition}\label{prop:topology}
	Let $\scr C$ be an $\infty$-category and let $\scr E$ be a collection of presheaves on $\scr C$. Let $\tau$ be the finest quasi-topology on $\scr C$ such that $\scr E\subset \Shv_\tau(\scr C)$. Then $\tau$ is a topology.
	\end{proposition}

\begin{proof}
	To begin with, note that $\tau$ exists: for $X\in\scr C$, $\tau(X)$ is the collection of sieves $R\into X$ such that, for every $f\colon Y\to X$ in $\scr C$ and every $F\in \scr E$, the map
	\[\Map(Y,F)\to\Map(f^\ast R,F)\]
	is an equivalence.
	To prove that $\tau$ is a topology, we must verify that, if $S\in\tau(X)$ and $R$ a sieve on $X$ such that $g^\ast R\in\tau(X')$ for every $g\colon X'\to X$ in $S$, then $R\in\tau(X)$. Let $f\colon Y\to X$ be a morphism in $\scr C$ and let $F\in \scr E$. We must show that the left vertical arrow in the square
	\begin{tikzmath}
		\def\colsep{1em}
		\diagram{\Map(Y,F) & \Map(f^\ast S,F) \\ \Map(f^\ast R,F) & \Map(f^\ast S\times_Yf^\ast R,F) \\};
		\arrows (11-) edge (-12) (11) edge (21) (21-) edge (-22) (12) edge (22);
	\end{tikzmath}
	is an equivalence. We will show that the other three arrows are equivalences. The top horizontal arrow is an equivalence because $S\in\tau(X)$. For the right vertical arrow, write $f^\ast S\simeq\colim_{Z\in\scr C/f^\ast S}Z$ as a (possibly large) colimit of representables. Since colimits in $\PSh(\scr C)$ are universal, $f^\ast R\times_Yf^\ast S\simeq \colim_{Z\in\scr C/f^\ast S} f^\ast R\times_Y Z$. For every $Z\to f^\ast S$, $f^\ast R\times_Y Z$ belongs to $\tau(Z)$ by assumption, and hence
	\[\Map(f^\ast S,F)\simeq\lim_Z\Map(Z,F)\simeq\lim_Z\Map(f^\ast R\times_Y Z,F)\simeq \Map(f^\ast R\times_Yf^\ast S,F).\]
	The proof that the bottom horizontal arrow is an equivalence is similar: write $f^\ast R\simeq\colim_{Z\in\scr C/f^\ast R}Z$ and use that $f^\ast S\in\tau(Y)$.
\end{proof}

\begin{corollary}\label{cor:gentop}
	Let $\scr C$ be an $\infty$-category and $\tau$ a quasi-topology on $\scr C$. Then $\Shv_\tau(\scr C)=\Shv_{\bar\tau}(\scr C)$.
	 \end{corollary}

\begin{proof}
	Note that $\Shv_{\bar\tau}(\scr C)\subset\Shv_\tau(\scr C)$. Let $\rho$ be the finest quasi-topology on $\scr C$ such that $\Shv_\tau(\scr C)\subset \Shv_\rho(\scr C)$. Tautologically, $\rho$ contains $\tau$. By Proposition~\ref{prop:topology}, $\rho$ contains $\bar\tau$. Hence, $\Shv_\rho(\scr C)\subset \Shv_{\bar\tau}(\scr C)$.
	\end{proof}

We will also need an easy-to-use version of the ``comparison lemma'' \cite[III, Théorème 4.1]{SGA4-1} for sheaves of spaces:

\begin{lemma}\label{lem:comparison}
	Let $\scr D$ be an $\infty$-category, $\scr C$ a \emph{small} $\infty$-category, and $u\colon\scr C\into \scr D$ a fully faithful functor. Let $\tau$ and $\rho$ be quasi-topologies on $\scr C$ and $\scr D$, respectively. Suppose that:
	\begin{enumerate}
		\item[(a)] Every $\tau$-sieve is generated by a cover $\{U_i\to X\}$ such that:
		\begin{enumerate}
			\item[(a1)] the fiber products $U_{i_0}\times_X\dotsb\times_XU_{i_n}$ exist and are preserved by $u$;
			\item[(a2)] $\{u(U_i)\to u(X)\}$ is a $\bar\rho$-cover $\scr D$.
		\end{enumerate}
		\item[(b)] For every $X\in\scr C$ and every $\rho$-sieve $R\into u(X)$, $u^*(R)\into X$ is a $\bar\tau$-sieve in $\scr C$.
		\item[(c)] Every $X\in\scr D$ admits a $\bar\rho$-cover $\{U_i\to X\}$ such that the fiber products $U_{i_0}\times_X\dotsb\times_X U_{i_n}$ exist and belong to the essential image of $u$.
	\end{enumerate}
	Then the adjunction $u^*\dashv u_*$ restricts to an equivalence of $\infty$-categories $\Shv_\rho(\scr D)\simeq\Shv_\tau(\scr C)$.
		\end{lemma}

We can rephrase the conclusion of the lemma as follows: a presheaf on $\scr D$ is a $\rho$-sheaf iff it is the right Kan extension of a $\tau$-sheaf on $\scr C$. An immediate consequence of the lemma is that the inclusion $\Shv_\rho(\scr D)\subset\PSh(\scr D)$ admits a left exact left adjoint $a_\rho$, namely the composition $u_*a_\tau u^*$.

\begin{proof}
	We tacitly use Corollary~\ref{cor:gentop} throughout the proof.
	We first show that $u^*$ and $u_*$ preserve sheaves. Let $\mathfrak U$ be a $\tau$-cover as in (a) and let \[\check C(\mathfrak U)\in\Fun(\Delta^\op,\PSh(\scr C))\] be its Čech nerve (note that $\colim\check C(\mathfrak U)$ is the sieve generated by $\mathfrak U$). By (a1), $u_!\check C(\mathfrak U)\simeq \check C(u(\mathfrak U))$, and by (a2), $u(\mathfrak U)$ is a $\bar\rho$-cover. 
	If $F$ is a $\rho$-sheaf, we deduce that
	\[\Map(u_!X, F)\to \Map (u_!\colim\check C(\mathfrak U), F)\]
		is an equivalence. By adjunction, $u^*$ preserves sheaves.
	Let $X\in\scr D$ and let $R\into X$ be a $\rho$-sieve. We claim that $u^*(R)\into u^*(X)$ becomes an equivalence after $\tau$-sheafification. By the universality of colimits in $\PSh(\scr C)$, it suffices to show that, for every $Y\in\scr C$ and every morphism $u(Y)\to X$, $u^*(R\times_Xu(Y)) \into Y$ is a $\bar\tau$-sieve. This follows from (b) since $R\times_Xu(Y)$ is a $\rho$-sieve. By adjunction, $u_*$ preserves sheaves. Thus, the adjunction $u^*\dashv u_*$ restricts to an adjunction
	\begin{tikzmath}
		\def\colsep{.9em}
		\diagram{u^*:\Shv_\rho(\scr D) & \Shv_\tau(\scr C):u_* \\};
		\arrows (11-) edge[vshift=\dbl] (-12)
		(-12) edge[c->,vshift=\dbl] (11-);
	\end{tikzmath}
	where $u_*$ is fully faithful. It remains to show that $u^*$ is conservative on $\Shv_\rho(\scr D)$, but this follows at once from (c).
	\end{proof}

A cartesian square of schemes
\begin{tikzmath}
	\diagram{W & V \\ U & X \\};
	\arrows (11-) edge[c->] (-12) (11) edge (21) (21-) edge[c->] node[above]{$j$} (-22) (12) edge node[right]{$p$} (22);
\end{tikzmath}
will be called a \emph{Nisnevich square} over $X$ if $j$ is an open immersion, $p$ is étale, and there exists a closed immersion $Z\into X$ complement to $U$ such that $p$ induces an isomorphism $V\times_XZ\simeq Z$.
We say that such a square is \emph{finitely presented} if $j$ and $p$ are finitely presented.

Let $B$ be a scheme. We denote by $\Sm_B$ the category of smooth $B$-schemes and by $\Sm_B'\subset \Sm_B$ the full subcategory spanned by compositions of open immersions and finitely presented smooth morphisms. If $B$ is coherent, we also consider the subcategory $\Sm_B^\fp\subset\Sm_B$ of finitely presented smooth $B$-schemes.
We will define the following quasi-topologies on $\Sm_B$: \begin{tikzmath}
	\def\colsep{1em}
	\def\rowsep{1em}
	\diagram{ & & \Nis_\qc^\fp \\
	\Zar & & \Nis_\qc \\
	& \Nis\rlap. & \\};
	\path (13) edge (23) (23) edge (32) (21) edge (32)
	;
\end{tikzmath}
The quasi-topology $\Zar$ will also be defined on $\Sm_B'$, and $\Nis_\qc^\fp$ and $\Nis$ will also be defined on $\Sm_B'$ and $\Sm_B^\fp$.
The $\Zar$-sieves are the sieves generated by open covers.
 The quasi-topology $\Nis_\qc$ (\resp{} $\Nis_\qc^\fp$) consists of:
\begin{itemize}
	\item the empty sieve on $\emptyset$;
	\item for every Nisnevich square (\resp{} finitely presented Nisnevich square) as above, the sieve generated by $\{j,p\}$.
\end{itemize}
The \emph{Nisnevich quasi-topology} $\Nis$ is then defined as follows on each category: $\Nis=\Zar\cup\Nis_\qc$ on $\Sm_B$, $\Nis=\Zar\cup\Nis_\qc^\fp$ on $\Sm_B'$, and $\Nis=\Nis_\qc^\fp$ on $\Sm_B^\fp$.

\begin{lemma}\label{lem:Nisfp}
	For every $\Nis_\qc$-sieve $R\into X$ in $\Sm_B$, there exists an open cover $\{f_i\colon X_i\into X\}$ such that $f_i^*R$ contains a $\Nis_\qc^\fp$-sieve.
	\end{lemma}

\begin{proof}
	Let $j:U\into X\from V:p$ be a Nisnevich square generating $R$, and let $Z$ be a closed complement of $j$ such that $V\times_XZ\simeq Z$.
	Taking an open cover of $X$ if necessary, we may assume that $X$ is coherent. Let $\{V_i\}$ be an open cover of $V$ by coherent schemes, and let $X_i=p(V_i)$. Then $V_i\to X_i$ is finitely presented and is an isomorphism over $Z\cap X_i$. Since $\{U,X_i\}$ is an open cover of $X$, we may assume that $p$ is finitely presented.
	As $X$ is coherent, we can write $Z=\lim_\alpha Z_\alpha$ where each $Z_\alpha$ is a finitely presented closed subscheme of $X$. Since $p$ is finitely presented and is an isomorphism over $Z$, there exists $\alpha$ such that $p$ is an isomorphism over $Z_\alpha$. If $j_\alpha$ is the open immersion complement to $Z_\alpha$, then $\{j_\alpha,p\}$ is a $\Nis_\qc^\fp$-cover refining $\{j,p\}$. 
	 \end{proof}

We say that a presheaf $F$ on $\Sm_B^{(\fp)}$ satisfies \emph{Nisnevich excision} if:
\begin{itemize}
	\item $F(\emptyset)\simeq *$;
	\item for every Nisnevich square $Q$ in $\Sm_B^{(\fp)}$, $F(Q)$ is cartesian.
\end{itemize}

\begin{proposition}\label{prop:Nistopos}
	Let $B$ be a scheme. \begin{enumerate}
		\item A presheaf on $\Sm_B$ is a Nisnevich sheaf if and only if it is the right Kan extension of a Nisnevich sheaf on $\Sm_B'$. In particular, $\Shv_\Nis(\Sm_B)$ is an $\infty$-topos and the inclusion $\Shv_\Nis(\Sm_B)\subset\PSh(\Sm_B)$ admits a left exact left adjoint.
		\item A presheaf on $\Sm_B$ is a Nisnevich sheaf if and only if it satisfies Zariski descent and Nisnevich excision.
	\end{enumerate}
	If $B$ is coherent, then:
	\begin{enumerate}
		\setcounter{enumi}{2}
		\item A presheaf on $\Sm_B$ is a Nisnevich sheaf if and only if it is the right Kan extension of a Nisnevich sheaf on $\Sm_B^\fp$. In particular, $\Shv_\Nis(\Sm_B)\simeq\Shv_\Nis(\Sm_B^\fp)$.
		\item A presheaf on $\Sm_B^\fp$ is a Nisnevich sheaf if and only if it satisfies Nisnevich excision.
	\end{enumerate}
\end{proposition}

\begin{proof}
	(1) It suffices to verify the assumptions of Lemma~\ref{lem:comparison} for the inclusion $\Sm_B'\subset\Sm_B$. The only nontrivial point is (b), which follows from Lemma~\ref{lem:Nisfp}. 
	(3) By (1), it suffices to verify the assumptions of Lemma~\ref{lem:comparison} for the inclusion $\Sm_B^\fp\subset\Sm_B'$. For (c), note that every scheme in $\Sm_B'$ is quasi-separated.

	(2,4) For $Q$ a Nisnevich square $j:U\into X\from V:p$ in $\Sm_B^{(\fp)}$, denote by $C_Q\in\PSh(\Sm_B^{(\fp)})$ the colimit of the Čech nerve $\check C(\{j,p\})$ (\ie, the sieve generated by $\{j,p\}$) and by $K_Q\in \PSh(\Sm_B^{(\fp)})$ the pushout of $Q$. Let $C$ (\resp{} $K$) be the class of morphisms of the form $C_Q\to X$ (\resp{} $K_Q\to X$) in $\PSh(\Sm_B^{(\fp)})$, where $Q$ is any Nisnevich square, together with the empty sieve on the empty scheme. By definition, a presheaf is a $\Nis_\qc^{(\fp)}$-sheaf iff it is $C$-local, and it satisfies Nisnevich excision iff it is $K$-local. 
	The arguments of \cite[\S5]{VV:cd} show that $C$ and $K$ generate the same class of morphisms under 2-out-of-3, pushouts, and colimits, whence the result.
	\end{proof}

These technical preliminaries aside, we can now define the unstable motivic homotopy category $\H(B)$ of an arbitrary scheme $B$. We say that a presheaf $F$ on $\Sm_B$ is \emph{$\A^1$-invariant} if, for every $X\in\Sm_B$, the projection $\A^1\times X\to X$ induces an equivalence $F(X)\simeq F(\A^1\times X)$. Note that if $F$ is Nisnevich sheaf on $\Sm_B$, the $\A^1$-invariance condition can be checked on $\Sm_B'$, and even on $\Sm_B^\fp$ if $B$ is coherent.

We let $\H(B)\subset\Shv_\Nis(\Sm_B)$ be the full subcategory of $\A^1$-invariant Nisnevich sheaves.
This definition is of course equivalent to the standard one when $B$ is noetherian and of finite Krull dimension.
By Proposition \ref{prop:Nistopos} (1), $\H(B)$ is a presentable $\infty$-category and the inclusion $\H(B)\subset \PSh(\Sm_B)$ admits a left adjoint
\[M\colon \PSh(\Sm_B)\to \H(B).\]

\begin{proposition}\label{prop:products}
	The functor $M$ preserves finite products.
\end{proposition}

\begin{proof}
	As $M$ factors through $\PSh(\Sm_B')$, it suffices to show that $M'\colon \PSh(\Sm_B')\to \H(B)$ preserves finite products. The functor
	\[L_{\A^1}\colon F\mapsto \colim_{n\in\Delta^\op} F(\A^n\times\ph)\]
	is left adjoint to the inclusion of $\A^1$-invariant presheaves into all presheaves, and it preserves finite products since $\Delta^\op$ is sifted. Let $a_\Nis$ be the Nisnevich sheafification functor. A standard argument shows that there exists an ordinal $\alpha$ such that the $\alpha$th iteration of $L_{\A^1}\circ a_\Nis$, viewed as a pointed endofunctor of $\PSh(\Sm_B')$, is equivalent to $M'$. Since $L_{\A^1}$, $a_\Nis$, and transfinite composition preserve finite products, so does $M'$.
\end{proof}

As is usual, if $X\in\Sm_B$, we will commit an abuse of notation and denote by $X$ the image of $X$ by the functor $\Sm_B\to\H(B)$, composition of the Yoneda embedding and the localization functor $M$. 
If $f\colon B'\to B$ is a morphism of schemes, the base change functor $\Sm_B\to \Sm_{B'}$ preserves trivial line bundles and Čech nerves of Nisnevich covers. It follows that the functor \[\PSh(\Sm_{B'})\to\PSh(\Sm_B),\quad F\mapsto F(\ph\times_BB'),\] preserves $\A^1$-invariant Nisnevich sheaves, and hence restricts to a limit-preserving functor $f_*\colon \H(B')\to \H(B)$. We denote by $f^*$ its left adjoint; it preserves finite products by Proposition~\ref{prop:products}. If $f$ is smooth, the base change functor $\Sm_B\to\Sm_{B'}$ has a left adjoint, namely the forgetful functor $\Sm_{B'}\to\Sm_B$, which also preserves trivial line bundles and Čech nerves of Nisnevich covers. It follows that in this case $f^*$ has a left adjoint $f_\sharp\colon\H(B')\to\H(B)$. We immediately verify that the exchange transformation $\Ex_\sharp^*$ and the projector $\Pr_\sharp^*$ are equivalences.

\begin{proposition}\label{prop:H(B)}
	\leavevmode
	\begin{enumerate}
		\item If $B$ is a coherent scheme, every $X\in\Sm_B^\fp$ is compact in $\H(B)$.
		\item If $f\colon B'\to B$ is coherent, $f_*\colon \H(B')\to \H(B)$ preserves filtered colimits.
		\item If $B$ is the limit of a cofiltered diagram of coherent schemes $B_\alpha$ with affine transition maps, then $\H(B)\simeq\lim_\alpha \H(B_\alpha)$ in the $\infty$-category of $\infty$-categories.
	\end{enumerate}
\end{proposition}

\begin{proof}
	(1) It suffices to show that $\H(B)$ is closed under filtered colimits in $\PSh(\Sm_B^\fp)$. In fact, it is obvious that the subcategories of $\A^1$-invariant presheaves and of presheaves satisfying Nisnevich excision are both closed under filtered colimits.
	
	(2) By Zariski descent, we can assume that $B$ is coherent. The $\infty$-category $\H(B)$ is then generated under colimits by $X\in\Sm_B^\fp$. By (1), these generators are compact and are carried by $f^*$ to compact objects in $\H(B')$. The result follows by adjunction.
	
	(3) In this situation, the category $\Sm_B^\fp$ is the colimit of the categories $\Sm_{B_\alpha}^\fp$, and hence $\PSh(\Sm_B^\fp)$ is the limit of the $\infty$-categories $\PSh(\Sm_{B_\alpha}^\fp)$. It remains to show that $F\in\PSh(\Sm_B^\fp)$ is an $\A^1$-invariant Nisnevich sheaf if, for all $\alpha$, its image in $\PSh(\Sm_{B_\alpha}^\fp)$ is. This follows from the fact that any trivial line bundle (\resp{} any Nisnevich square) in $\Sm_B^\fp$ is the pullback of a trivial line bundle (\resp{} a Nisnevich square) in $\Sm_{B_\alpha}^\fp$, for some $\alpha$.
\end{proof}

Our next goal is to generalize the gluing theorem of Morel–Voevodsky \cite[Theorem 3.2.21]{MV} to our setting. The proof in \textit{loc.\ cit.} uses the fact that henselian local schemes form a conservative family of points of the $\infty$-topos $\Shv_\Nis(\Sm_B)$, which is not true anymore when $B$ is not noetherian of finite Krull dimension. While it is not difficult to give a proof of the gluing theorem that avoids the use of points and works in general, we will give instead a shorter argument that reduces the general case to the Morel–Voevodsky case.

Suppose that $B$ is a cofiltered limit of coherent schemes $B_\alpha$, so that $\H(B)\simeq\lim_\alpha\H(B_\alpha)$. Let $f_{\beta\alpha}\colon B_\beta\to B_\alpha$ be the transition maps and $f_\alpha\colon B\to B_\alpha$ the canonical projections. Then, by \cite[Lemma 6.3.3.6]{HTT},
\begin{equation}\label{eqn:continuous1}
	\id_{\H(B)}\simeq \colim_\alpha f_\alpha^*f_{\alpha *}.
\end{equation}
Moreover, since functors of the form $f_{*}$ preserve filtered colimits, the left adjoint functors $f_\alpha^*$ can be computed as follows:
\begin{equation}\label{eqn:continuous2}
f_{\beta*}f_\alpha^*\simeq\colim_\gamma f_{\gamma\beta*}f_{\gamma\alpha}^*.
\end{equation}

\begin{samepage}
\begin{proposition}
	\label{prop:locality}
	Let $B$ be a scheme and let $i\colon Z\into B$ be a closed immersion with open complement $j\colon U\into B$. Then:
	\begin{enumerate}
	\item For every $F\in\H(B)$, the square
	\begin{tikzmath}
		\diagram{j_\sharp j^*F & F \\ j_\sharp j^*B & i_*i^* F \\};
		\arrows (11-) edge node[above]{$\epsilon$} (-12) (11) edge (21) (21-) edge (-22) (12) edge node[right]{$\eta$} (22);
	\end{tikzmath}
	is cocartesian.
	\item $i_*\colon\H(Z)\to\H(B)$ is fully faithful. \end{enumerate}
\end{proposition}
\end{samepage}

\begin{proof}
	(1) By Zariski descent, we can assume that $B$ is coherent. Let $\{i_\alpha\colon Z_\alpha\into B\}$ be the cofiltered poset of finitely presented closed subschemes of $B$ containing $Z$, and let $j_\alpha\colon U_\alpha\into B$ be the open immersion complement to $i_\alpha$. Then $Z\simeq\lim_\alpha Z_\alpha$ and $\{U_\alpha\}$ is an open cover of $U$ which is closed under binary intersections. By Zariski descent, the canonical transformation
	\[\colim_\alpha j_{\alpha\sharp} j_\alpha^*\to j_\sharp j^*\]
	is an equivalence. On the other hand, by~\eqref{eqn:continuous2} and a cofinality argument, the canonical transformation
	\[\colim_\alpha i_{\alpha*}i_\alpha^*\to i_*i^*\]
	is an equivalence. It therefore suffices to prove the result when $B$ is coherent and $i$ is finitely presented.
	In that case, we can write $i$ and $j$ as cofiltered limits of complementary immersions $i_\alpha\colon Z_\alpha\into B_\alpha$ and $j_\alpha\colon U_\alpha\into B_\alpha$,
	such that $B_\alpha$ is of finite type over $\Z$ and such that the squares
		\begin{tikzmath}
			\diagram{Z  & B & U \\ Z_\alpha & B_\alpha & U_\alpha \\};
			\arrows (11-) edge[c->] node[above]{$i$} (-12) (21-) edge[c->] node[below]{$i_\alpha$} (-22)
			(-13) edge[left hook->] node[above]{$j$} (12-) (-23) edge[left hook->] node[below]{$j_\alpha$} (22-)
			(11) edge node[left]{$f_\alpha$} (21) (12) edge node[left]{$f_\alpha$} (22) (13) edge node[left]{$f_\alpha$} (23);
		\end{tikzmath}
		are cartesian. Let $F_\alpha=f_{\alpha*}F$ be the component of $F$ in $\H(B_\alpha)$, so that, by~\eqref{eqn:continuous1}, $F\simeq\colim_\alpha f_\alpha^*F_\alpha$. Since $B_\alpha$ is noetherian of finite Krull dimension, we have a cocartesian square
		\begin{tikzmath}
			\diagram{j_{\alpha\sharp} j_\alpha^*F_\alpha & F_\alpha \\ j_{\alpha\sharp} U_\alpha & i_{\alpha*}i_\alpha^* F_\alpha \\};
			\arrows (11-) edge node[above]{$\epsilon$} (-12) (11) edge (21) (21-) edge (-22) (12) edge node[right]{$\eta$} (22);
		\end{tikzmath}
		in $\H(B_\alpha)$. Applying $f_\alpha^*$ and taking the colimit over $\alpha$, we obtain a cocartesian square in $\H(B)$ which maps canonically to the given square. Moreover, the maps on the top left, bottom left, and top right corners are equivalences since $f_\alpha^*j_{\alpha\sharp}\simeq j_\sharp f_\alpha^*$ and since $j_\sharp$ and $j^*$ preserve colimits. It remains to prove that the map
		\[\colim_\alpha f_\alpha^*i_{\alpha*}i_\alpha^*F_\alpha\to i_*i^*F\]
		on the bottom right corner	is an equivalence. Since $i_{*}$ preserves filtered colimits, it suffices to show that the exchange transformation $\Ex^*_*\colon f_\alpha^*i_{\alpha*}\to i_*f_\alpha^*$ is an equivalence. Using~\eqref{eqn:continuous2} and the fact that the exchange transformation $f_{\gamma\alpha}^*i_{\alpha*}\to i_{\gamma*}f_{\gamma\alpha}^*$ is an equivalence (which is a consequence of the gluing theorem for finite-dimensional noetherian schemes), we compute:
		\begin{multline*}
			f_{\beta*}f_\alpha^*i_{\alpha*}\simeq\colim_{\gamma}f_{\gamma\beta*}f_{\gamma\alpha}^*i_{\alpha*}\simeq \colim_\gamma f_{\gamma\beta*} i_{\gamma*}f_{\gamma\alpha}^*\\
			\simeq \colim_\gamma i_{\beta*}f_{\gamma\beta*}f_{\gamma\alpha}^*\simeq i_{\beta*}\colim_\gamma f_{\gamma\beta*}f_{\gamma\alpha}^*\simeq i_{\beta*}f_{\beta*}f_\alpha^*\simeq f_{\beta*}i_*f_\alpha^*.
		\end{multline*}
		One verifies easily that this composition coincides with $f_{\beta*}\Ex^*_*$, which completes the proof.
	
		(2) Applying (1) to $i_*F$, we deduce that the unit $\id\to i_*i^*$ is an equivalence on $i_*F$. It follows from a triangle identity that the counit $i^*i_*\to \id$ becomes an equivalence after applying $i_*$. By \cite[Proposition 18.1.1]{EGA4-4}, $\H(Z)$ is generated under colimits by pullbacks of smooth $B$-schemes. It follows that $i_*$ is conservative and hence fully faithful. \end{proof}

Denote by $\H_*(B)$ the undercategory $\H(B)_{B/}$. All the features of $\H(B)$ discussed so far have obvious analogs for $\H_*(B)$. The smash product $\wedge$ on $\H_*(B)$ is the unique symmetric monoidal product which is compatible with colimits and for which the functor $(\ph)_+\colon \H(B)\to\H_*(B)$ is symmetric monoidal.
One can then define the $\infty$-category $\SH(B)$ as a symmetric monoidal presentable $\infty$-category as in \cite[Definition 4.8]{Robalo}, by formally inverting $S^{\A^1}$ for the smash product on $\H_*(B)$. We thus have a symmetric monoidal colimit-preserving functor
 \[\Sigma^\infty\colon\H_*(B)\to\SH(B),\]
and we let $\Sigma^\infty_+=\Sigma^\infty\circ(\ph)_+$.
Note that $\SH(B)$ is stable since $S^{\A^1}$ is the suspension of the pointed motivic space $(\A^1\minus 0,1)$. Because the cyclic permutation of $S^{\A^1}\wedge S^{\A^1}\wedge S^{\A^1}$ is homotopic to the identity, $\SH(B)$ can also be described as the following limit of $\infty$-categories:
\begin{equation}\label{eqn:defSH}
	\SH(B)=\lim(\dotsb\xrightarrow{\Omega^{\A^1}}\H_*(B)\xrightarrow{\Omega^{\A^1}}\H_*(B)),
\end{equation}
where $\Omega^{\A^1}$ is right adjoint to $\Sigma^{\A^1}$ \cite[Corollary 4.24]{Robalo}.

If $f\colon B'\to B$ is a morphism of schemes, then $f_*\Omega^{\A^1}\simeq\Omega^{\A^1}f_*$ and hence $f_*$ induces a limit-preserving functor $f_*\colon\SH(B')\to\SH(B)$. Its left adjoint $f^*$ is the unique colimit-preserving symmetric monoidal functor $f^*\colon \SH(B)\to \SH(B')$ such that $f^*\Sigma^\infty_+ X=\Sigma^\infty_+(X\times_BB')$ for $X\in\Sm_B$.

\begin{proposition}\label{prop:SH(B)}
	\leavevmode
	\begin{enumerate}
		\item $\SH(B)$ is generated under colimits by objects of the form $\Sigma^{-\A^n}\Sigma^\infty_+X$ for $X\in\Sm_B$ and $n\geq 0$.
		\item If $B$ is a coherent scheme and $X\in\Sm_B^\fp$, $\Sigma^\infty_+X\in\SH(B)$ is compact.
		\item If $f\colon B'\to B$ is coherent, $f_*\colon \SH(B')\to \SH(B)$ preserves colimits (and hence admits a right adjoint).
		\item If $B$ is the limit of a cofiltered diagram of coherent schemes $B_\alpha$ with affine transition maps, then $\SH(B)\simeq\lim_\alpha \SH(B_\alpha)$ in the $\infty$-category of $\infty$-categories.
	\end{enumerate}
\end{proposition}

\begin{proof}
	(1) Let $E\in\SH(B)$ have components $E_n\in \H_*(B)$. By \cite[Lemma 6.3.3.6]{HTT},
	\[E\simeq\colim_{n\geq 0}\Sigma^{-\A^n}\Sigma^\infty E_n,\]
	and each $\Sigma^{-\A^n}\Sigma^\infty E_n$ is clearly an iterated colimit of objects of the desired form.
	
	(2) By Proposition~\ref{prop:H(B)} (1), $\H_*(B)$ is compactly generated by $X_+$, $X\in\Sm_B^\fp$. The object $S^{\A^1}\in \H_*(B)$ is compact, being a finite colimit of compact objects, and so the functor $\Omega^{\A^1}\colon \H_*(B)\to\H_*(B)$ preserves filtered colimits. The assertion now follows immediately from~\eqref{eqn:defSH}.
	
	(3) We can assume that $B$ is coherent. By (1) and (2), $f^*$ sends a family of compact generators of $\SH(B)$ to compact objects in $\SH(B')$. By adjunction, $f_*$ preserves filtered colimits. Since $f_*$ preserves limits and both $\SH(B')$ and $\SH(B)$ are stable, it also preserves finite colimits.
	
	(4) This follows from Proposition~\ref{prop:H(B)} (3) and~\eqref{eqn:defSH}.
\end{proof}

Finally, we prove that $\SH(\ph)$ satisfies the proper base change theorem and related properties:

\begin{proposition}
	\label{prop:PBC}
	Let
	\begin{tikzmath}
		\diagram{Y' & Y \\ X' & X \\};
		\arrows (11-) edge node[above]{$g$} (-12) (12) edge node[right]{$p$} (22) (11) edge node[left]{$p'$} (21) (21-) edge node[below]{$f$} (-22);
	\end{tikzmath}
	be a cartesian square of  schemes where $p$ is proper.
	\begin{enumerate}
		\item For every $E\in\SH(Y)$, the exchange transformation $\Ex_*^*\colon f^*p_* E\to p'_*g^* E$ is an equivalence.
		\item For every $E\in\SH(Y)$ and $F\in\SH(X)$, the projector $\Pr_*^*\colon p_* E\wedge F\to p_*(E\wedge p^*F)$ is an equivalence.
		\item Suppose that $f$ is smooth. For every $E\in\SH(Y')$, the exchange transformation $\Ex_{\sharp *}\colon f_\sharp p'_* E\to p_* g_\sharp E$ is an equivalence.
	\end{enumerate}
\end{proposition}

\begin{proof}
	If $p$ is a closed immersion, all three statements follow easily from the gluing theorem. The argument of \cite[\S1.7.2]{Ayoub} shows that
	the map $p_\sharp\to p_*\pi_{1\sharp}\delta_*$ induced by $\Ex_{\sharp*}$
	is an equivalence when $p$ is a projection $\P^n_X\to X$.
	The proof of \cite[Lemma 2.4.23]{CD} then shows that (1–3) hold for such $p$. By Zariski descent, one immediately deduces (1–3) for $p$ projective. It remains to extend the results to $p$ proper. 
	(1) By Zariski descent, we can assume that $X$ and $X'$ are coherent. Let $C\colon \SH(Y)\to \SH(X')$ be the cofiber of the transformation $\Ex_*^*$. Since $\SH(X')$ is stable and compactly generated, it will suffice to show that $[K,C(E)]=0$ for every $E\in \SH(Y)$ and every $K\in\SH(X')$ compact. Fix $x\colon K\to C(E)$ and consider the poset $\Phi$ of closed subschemes $i\colon Z\into Y$ such that the image of $x$ in $[K,C(i_*i^*E)]$ is \emph{not} zero. If $\{i_\alpha\colon Z_\alpha\into Y\}$ is a cofiltered diagram of closed subschemes of $Y$ with limit $i\colon Z\into Y$, it follows from~\eqref{eqn:continuous2} that $\colim_\alpha i_{\alpha*}i_\alpha^* \simeq i_*i^*$. Since the source and target of $\Ex_*^*$ preserve filtered colimits, the canonical map
	\[\colim_\alpha C(i_{\alpha*}i_\alpha^*E)\to C(i_*i^*E)\]
	is an equivalence. By compactness of $K$, we deduce that $\Phi$ is closed under cofiltered intersections. On the other hand, using Chow's lemma \cite[XII, \S7]{SGA4-3}, the gluing theorem, and (1,3) for $p$ projective, we easily verify that $\Phi$ does not have a minimal element. Hence, $\Phi$ is empty.
	
	(2) Same proof as (1).

	(3) Arguing as in (1) proves the result when $f$ is coherent. It also proves that $\Ex_{\sharp *}g^*$ is an equivalence, whence the result when $f$ is an open immersion. Without loss of generality, assume now that $X$ is coherent.
	Then $\SH(Y')$ is generated under colimits by the images of $h_\sharp$ where $h$ is the pullback of the inclusion of an open subscheme of $X'$ which is coherent over $X$, so the general case follows. \end{proof}

By Nagata's compactification theorem \cite{Conrad} and Proposition~\ref{prop:PBC} (3), we can apply Deligne's gluing theory and define the exceptional adjunction
\[f_!: \SH(X)\rightleftarrows\SH(Y): f^!\]
at the level of triangulated categories, for $f\colon X\to Y$ a separated morphism of finite type between coherent schemes. Following \cite[\S2]{CD}, we then obtain the complete formalism of six operations for coherent schemes as described in \S\ref{sec:review}.

\begin{remark}\label{rmk:gluing}
	It is possible to define $\SH(\ph)$ as a contravariant functor from the category of schemes to the $\infty$-category of symmetric monoidal presentable $\infty$-categories. Using the $\infty$-categorical generalization of Deligne's gluing theory developed in \cite{LZ}, 
	one can define the exceptional adjunction $f_!\dashv f^!$, the natural transformation $f_!\to f_*$, and all the exchange transformations and projectors involving exceptional functors, at the level of $\infty$-categories (for $f$ a separated morphism of finite type between coherent schemes).
	Since $\SH(\ph)$ is a Zariski sheaf, one can
	 further define all this data for any morphism $f$ which is locally of finite type.
	Once this is done, the conventions set at the end of~\S\ref{sec:introduction} can be ignored altogether, and ``separated of finite type'' can be replaced everywhere by ``locally of finite type''.
\end{remark}